%% file: QC_ResStress_P2.tex
\documentclass[12pt]{amsart}

\usepackage{graphicx}
\usepackage{color}
\usepackage{xcolor}
\usepackage{amsfonts, amsmath, amssymb, amsbsy}
\usepackage{environments}
\usepackage{algorithm}
\usepackage{mathrsfs} 
\usepackage{hyperref}
\usepackage{ulem}
\usepackage[top=2.5cm,bottom=2.0cm,left=2.0cm,right=2.0cm]{geometry}
\newtheorem{assumption}{Assumption}[section]

\renewcommand{\thetheorem}{\arabic{section}.\arabic{theorem}}

\definecolor{yscol}{HTML}{6622AA}

\begin{document}
\title[Adaptive Blended A/C Coupling with P2-FEM]{A Posteriori Analysis and Adaptive Algorithms for \\ Blended Type Atomistic-to-Continuum Coupling with \\ Higher-Order Finite Elements}

\author{Yangshuai Wang }
\address{University of British Columbia, 1984 Mathematics Road, Vancouver, BC, Canada.}
\email{yswang2021@math.ubc.ca}

\date{\today}

\begin{abstract}
The efficient and accurate simulation of material systems with defects using atomistic-to-continuum (a/c) coupling methods is a topic of considerable interest in the field of computational materials science. To achieve the desired balance between accuracy and computational efficiency, the use of {\it a posteriori} analysis and adaptive algorithms is critical. In this work, we present a rigorous {\it a posteriori} error analysis for three typical blended a/c coupling methods: the blended energy-based quasi-continuum (BQCE) method, the blended force-based quasi-continuum (BQCF) method, and the atomistic/continuum blending with ghost force correction (BGFC) method. We employ first and second-order finite element methods (and potentially higher-order methods) to discretize the Cauchy-Born model in the continuum region. The resulting error estimator provides both an upper bound on the true approximation error and a lower bound up to a theory-based truncation indicator, ensuring its reliability and efficiency. Moreover, we propose an {\it a posteriori} analysis for the energy error. We have designed and implemented a corresponding adaptive mesh refinement algorithm for two typical examples of crystalline defects. In both numerical experiments, we observe optimal convergence rates with respect to degrees of freedom when compared to {\it a priori} error estimates.

\end{abstract}

\maketitle



\input{notation.tex}
\def\Use{\Us^{1,2}}
\def\M{\mathcal{M}}

\section{Introduction}
\label{sec:intro}

Atomistic-to-continuum (a/c) coupling methods are a type of concurrent multi-scale scheme \cite{TadmorMiller:2012, van2020roadmap} that have received significant attention from both the engineering and mathematical communities over the past two decades \cite{Ortiz:1995a, Miller:2008, 2009_SP_LC_BD_PB_CMAME, LuOr:acta, olson2014optimization, tembhekar2017automatic}. The fundamental principle behind a/c coupling methods involves applying a more accurate atomistic model in a localized defected region, while introducing the continuum model (such as the Cauchy-Born rule) for areas far from the defect cores. 

A comprehensive overview of various a/c coupling methods for defect simulation can be found in \cite{van2020roadmap, Miller:2008}, while \cite{LuOr:acta} provides a rigorous analysis of these methods. In this work, we mainly focus on the blended types of a/c coupling methods, including the BQCE (blended energy-based quasi-continuum) method \cite{2012-CMAME-optbqce, 2014-bqce}, the BQCF (blended force-based quasi-continuum) method \cite{MiLu:2011, Dobson:2008a}, and the recently developed BGFC (atomistic/continuum blending with ghost force correction) method \cite{OrZh:2016, fang2020blended}.

One of the fundamental challenges in a/c coupling methods is to determine the optimal assignment of atomistic and continuum regions, as well as the mesh structure, in order to achieve a (quasi-)optimal balance between accuracy and efficiency. {\it A priori} choices, although feasible, often result in sub-optimal distribution of computational resources and are only suitable for simple setups such as single point defects. Therefore, {\it a posteriori} analysis and corresponding adaptive algorithms are essential for the efficient implementation and simulation of a/c coupling methods for crystalline defects. However, the development of {\it a posteriori} analysis and adaptive algorithms for a/c coupling methods in two or three dimensions is relatively recent, and it serves as a cornerstone for the application of these methods in real-world material systems.

The first approach to a/c coupling methods from a heuristic engineered standpoint was the goal-oriented {\it a posteriori} error estimate. Various applications of this approach on different model problems can be found in \cite{2003_TO_SP_JCP, 2006_SP_PB_TO_IJMCE, 2006_TO_SP_AR_PM_SISC, 2007_MA_ML_IJMCEpdf, 2008_MA_ML_CMAME}. However, the reliability of the error estimator cannot be guaranteed. On the other hand, the residual based {\it a posteriori} error estimate provides a quantitative estimate of the error in a certain energy-norm. This approach was first analyzed in \cite{Ortner:qnl.1d, OrtnerWang:2014} for consistent a/c coupling methods in one dimension and extended to two dimensions for the GRAC method with nearest neighbor interaction in \cite{wang2018posteriori}, and then to the case of finite range interactions in \cite{liao2018posteriori}. However, these studies were limited by the computational cost incurred in evaluating the modeling residual. To address this issue, the authors of \cite{wang2022adapac1} proposed a theory-based approximation for the residual-based {\it a posteriori} error estimator that significantly reduces the computational effort. Although the error estimators presented in these studies have been shown to be reliable, their derivation and construction are considered ad hoc and complicated, as discussed in \cite{wang2018analysis}. Recently, a residual force-based {\it a posteriori} error estimator was developed for blended a/c coupling methods (BQCE, BQCF, and BGFC), which is relatively intuitive and has been shown to be reliable \cite{wang2022adapac2}. However, it is important to note that the aforementioned studies utilized $\mathcal{P}_1$ finite element to discretize the Cauchy-Born approximation, and the adaptive computations with higher-order finite element and the {\it a posteriori} analysis for energy error have not been explored in the literature to the best of our knowledge.

The aim of this study is to present an improved residual-based {\it a posteriori} error estimation method for blended atomistic/continuum (a/c) coupling methods with higher order finite elements. The proposed error estimator is developed by solving an auxiliary Poisson equation, where the interpolated residual forces serve as the source term. This construction is based on the Rieze representation of the residual in dual norm and builds upon our prior work on quantum mechanics/molecular mechanics (QM/MM) coupling methods~\cite{wang2020posteriori}. The resulting error estimator provides both upper and lower bounds for the true approximation error, subject to a controllable truncation indicator. An adaptive algorithm is developed that allows for real-time adjustments to the atomistic, blending, and continuous mesh regions, enabling the optimal balance between these regions to be achieved. We also examine the truncation error and define a theory-based indicator to control the computational domain size. Finally, an error estimator for energy error is presented. Our tests include different types of crystalline defects that have not been previously considered in the literature on adaptive a/c coupling methods. The numerical results demonstrate that our adaptive algorithm achieves optimal convergence rates consistent with the {\it a priori} error estimate.

The scope of this paper is primarily focused on three specific blended a/c coupling methods for Bravais lattices with point defects in two dimensions, with the aim of presenting clear ideas and adaptive strategies. However, the analysis and methodology proposed can be extended to other consistent multi-scale coupling methods, multilattice crystals, and problems in three dimensions. The author plan to investigate the generalization of these approaches in future research and provide a discussion of potential future developments in Section~\ref{sec: conclusion}.

\subsection{Outline}
\label{sec:sub:outl}

This paper is structured as follows: Section~\ref{sec:A-AC_model} introduces the atomistic and continuum models, along with three typical blended atomistic/continuum (a/c) coupling schemes, including the BQCE, BQCF, and BGFC methods, which are the focus of this work. In Section~\ref{sec:a_post}, we review the {\it a priori} error estimate for the blended a/c coupling methods with higher order finite elements, which has been previously established in the literature \cite{2014-bqce, colz2016}, providing a solid foundation for our analysis. We then present the rigorous residual-based {\it a posteriori} error estimates (Theorem~\ref{th:mainresult}), which offer both upper and lower bounds on the true approximation error, subject to a controlled truncation indicator. Section~\ref{sec: numexp} presents our approach to assigning local error contributions and proposing an adaptive algorithm for the blended a/c coupling methods, based on the theoretical results established in the previous section. We report on numerical results obtained using this adaptive algorithm for crystalline defects that have not been previously considered in the literature, and provide a thorough discussion and explanation. Our conclusions and suggestions for future research are presented in Section~\ref{sec: conclusion}. Additional supporting materials are provided in the Appendices.

\subsection{Notations}
\label{sec:sub:not}

We use the symbol $\langle\cdot,\cdot\rangle$ to denote an abstract duality
pairing between a Banach space and its dual space. The symbol $|\cdot|$ normally
denotes the Euclidean or Frobenius norm, while $\|\cdot\|$ denotes an operator
norm.
For the sake of brevity of notation, we will denote $A\backslash\{a\}$ by
$A\backslash a$, and $\{b-a~\vert ~b\in A\}$ by $A-a$.
For $E \in C^2(X)$, the first and second variations are denoted by
$\<\delta E(u), v\>$ and $\<\delta^2 E(u) v, w\>$ for $u,v,w\in X$.
For a finite set $A$, we will use $\#A$ to denote the cardinality of $A$.
For second order tensors $\mA$ and $\mB$, we denote $\mA:\mB=\sum_{i,j}\mA_{ij}\mB_{ij}$.
The closed ball with radius $r$ and center $x$ is denoted by $B_r(x)$, or $B_r$ if the center is the origin.
To further simplify notation we will often write $\lesssim$ to mean $\leq C$ as well as $\eqsim$ to mean both $\lesssim$ and $\gtrsim$.
We use the standard definitions and notations $L^p$, $W^{k,p}$, $H^k$ for Lebesgue and Sobolev spaces. In addition we define the homogeneous Sobolev spaces
$
	\dot{H}^k(\Omega) := \big\{ f \in H^k_{\rm loc}(\Omega) \,|\,
										\nabla^k f \in L^2(\Omega) \big\}.
$

\section{Atomistic-to-Continuum (A/C) Coupling Model}
\label{sec:A-AC_model}

In this section, we present the atomistic model and its three blended a/c coupling approximations, namely the energy-based blended quasi-continuum (BQCE) method, the force-based blended quais-continuum (BQCF) method and the blended ghost force correction (BGFC) method. For the sake of simplicity, we mainly consider the point defects in two dimensions. The extension to three dimensions and more complex defects such as dislocations and cracks will require substantial additional efforts~\cite{2017-bcscrew, 2019-crackbif}, which will be explored in our future work. We keep the presentation of this section as concise as possible since much of the details can be found in various earlier works \cite{LuOr:acta, 2014-bqce, colz2016, 2012-CMAME-optbqce, fang2020blended}.

\subsection{Atomistic model}
\label{sec:sub:atom}

\def\L{\Lambda}
\def\Nhd{\mathcal{N}}
\def\Rg{\mathcal{R}}
\def\mF{\textsf{F}}
\def\mA{\textsf{A}}
\def\R{\mathbb{R}}
\def\Z{\mathbb{Z}}
\def\Lhom{\Lambda^{\rm hom}}
\def\E{\mathcal{E}}
\def\T{\mathcal{T}}
\def\Us{\mathscr{U}}
\def\UsH{\mathscr{U}^{1,2}}
\def\a{{\rm a}}
\def\b{{\rm b}}
\def\c{{\rm c}}
\def\h{{\rm h}}
\def\Rcore{R_{\rm DEF}}

Let $\Lhom=\mA\Z^2$ with some non-singular matrix $\mA\in\R^{2\times 2}$ be a perfect single lattice possessing no defects and $\L\subset \R^2$ be the corresponding single lattice with some point defects.
Without loss of generality, for the case of a single defect, we assume that it is contained within a ball $B_{\Rcore}$ for $\Rcore>0$; that is, $\L \cap B_{\Rcore}$ is finite and $\L \setminus B_{\Rcore} = \Lhom \setminus B_{\Rcore}$. The case with multiple local defects can be similarly defined~\cite{wang2022adapac1}.  

The deformed configuration of the infinite lattice $\L$ is a map $y \in \Us$ with $\mathscr{U}:=\{v:\L\to \mathbb{R}^2\}$, and it can be decomposed as
\begin{eqnarray}\label{y-u}
y(\ell) = \ell + u(\ell),  \qquad\forall~\ell\in\Lambda.
\end{eqnarray}

We define the interaction neighborhood for each $\ell\in \L$ by $\Nhd_{\ell} := \{ \ell' \in \L~|~0<|\ell'-\ell| \leq r_{\rm cut} \}$ with a given cut-off radius $r_{\rm cut}$. We also denote the interaction range $\Rg_\ell := \{\ell'-\ell~|~\ell'\in \Nhd_\ell\}$. For each atom $\ell\in \L$, we define the finite difference stencil for $u:\L\to \mathbb{R}^d$
\begin{align*}
Du(\ell):= \{D_\rho u(\ell)\}_{\rho \in \Rg_\ell} :=\{u(\ell+\rho)-u(\ell)\}_{\rho \in \Rg_\ell}.
\end{align*} 

For $u\in\Us$, the nodal interpolant of $u$ with respect to the background mesh is denoted as $I_{\rm a} u$~\cite{colz2016, 2013-defects, 2014-bqce}. Identifying $u=I_{\rm a}u$, we can define the piecewise constant gradient $\nabla u = \nabla I_{\rm a} u : \R^d \rightarrow \R^{d \times d}$, and
therefore the corresponding functional space of finite-energy displacements is defined by
\begin{align}\label{space:UsH}
\UsH(\L) := \big\{u:\L\rightarrow\mathbb{R}^{d} ~\big\lvert~ \|\nabla u\|_{L^2(\R^d)} < \infty \big\}.
\end{align}
We also define the following subspace of compact displacements
\begin{align}\label{space:Uc}
\Us^{\rm c}(\L) := \big\{u:\L\rightarrow\mathbb{R}^{d} ~\big\lvert~ \exists~R >0~{\rm s.t.}~u = {\rm const}~{\rm in}~\L\setminus B_{R}\big\}.
\end{align}

The site potential is a collection of mappings $V_{\ell}:(\R^3)^{\Rg_\ell}\rightarrow\R$, which represents the energy distributed to each atomic site in $\L$. We refer to \cite[Section 2]{2013-defects} for the detailed discussions on the assumptions of general site potentials. In this work, we will use the well-known EAM (Embedded Atom Method) model \cite{Daw:1984a} throughout the numerical experiments (cf.~Section~\ref{sec: numexp}). 

Furthermore, we assume that $V_{\ell}({\bf 0})=0$ for all $\ell \in \L$. As a matter of fact, if $V({\bf 0})\neq 0$, then it can be replaced with $V_{\ell}(Du)\equiv V_{\ell}(Du) - V_{\ell}({\bf 0})$; that is, $V$ should be understood as a site energy difference. We can formally define the energy-difference functional
\begin{align}\label{energy-difference}
\E^{\a}(u) = \sum_{\ell\in\Lambda}V_{\ell}\big(Du(\ell)\big).
\end{align}
It was shown in \cite{2013-defects} that $\E(u)$ is well-defined on the space $\UsH(\L)$. 

We can now rigorously formulate the variational problem for the equilibrium state as
\begin{equation}\label{eq:variational-problem}
u^{\a} \in \arg\min \big\{ \E^{\rm a}(u),~u \in \Us^{1,2}(\L) \big\},
\end{equation}
where ``$\arg\min$'' is understood as the set of local minimizers. 

In order to conduct an analysis of the approximation error, a stronger stability condition is necessary. However, this condition can only be rigorously justified for certain simple physical models. Therefore, we present it here as an assumption.
\begin{assumption}\label{as:LS}
	$\exists~\gamma>0$, such that $
		\big\< \delta^2\E^{\rm a}(u^{\rm a}) v , v\big\> \geq \gamma \|v\|^2_{\UsH} \qquad\forall~v\in\UsH(\L).
$
\end{assumption}

\subsection{A/C coupling methods}
\label{sec:sub:ac}

The atomistic-to-continuum (a/c) coupling methods is a class of concurrent multiscale methods which hybridizes atomistic and continuum models, aiming to achieve an optimal balance between accuracy and computational complexity. We refer to~\cite{van2020roadmap, luskin2013atomistic, miller2009unified} for the extensive overview and benchmark. In this work, we mainly focus on the blended types of a/c coupling methods since they are of practical interest and meanwhile outperform most alternative approaches~\cite{colz2016}.

\subsubsection{Domain decomposition and solution space}
\label{sec:sub:sub:dom}

Due to the infinite nature of the space in which \eqref{eq:variational-problem} is defined, it is impractical to solve it computationally. A common approach to obtain a feasible approximation is to impose a clamped boundary condition. Therefore, we limit the infinite lattice $\Lambda$ to a finite computational domain $\Omega$ with a radius of $R_{\Omega}$.

\begin{figure}[htb]
\begin{center}
	\includegraphics[height=6cm]{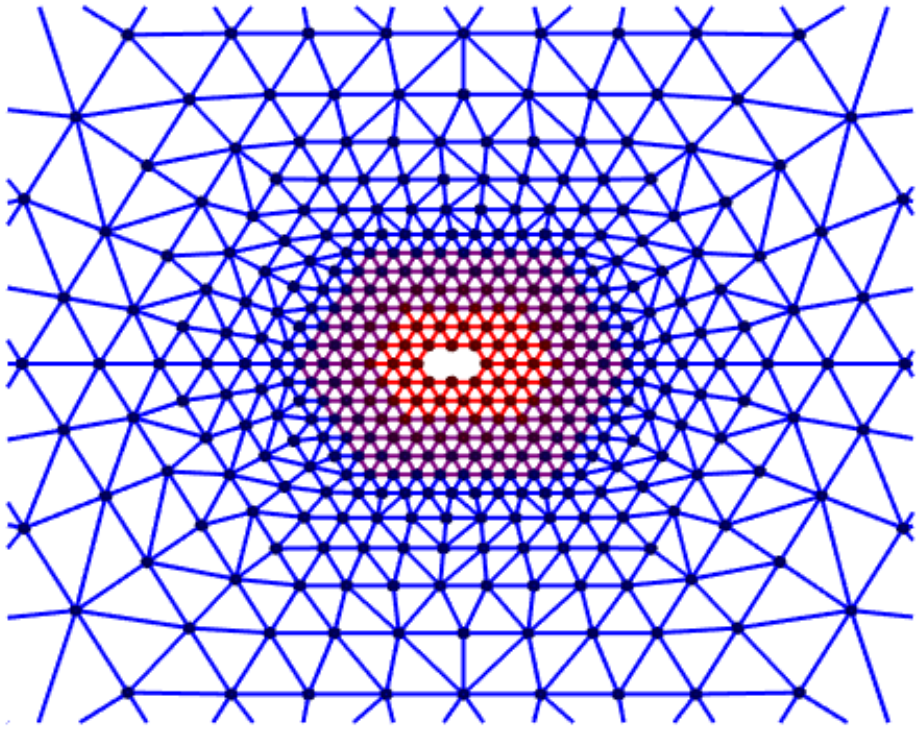}\quad
 \includegraphics[height=6cm]{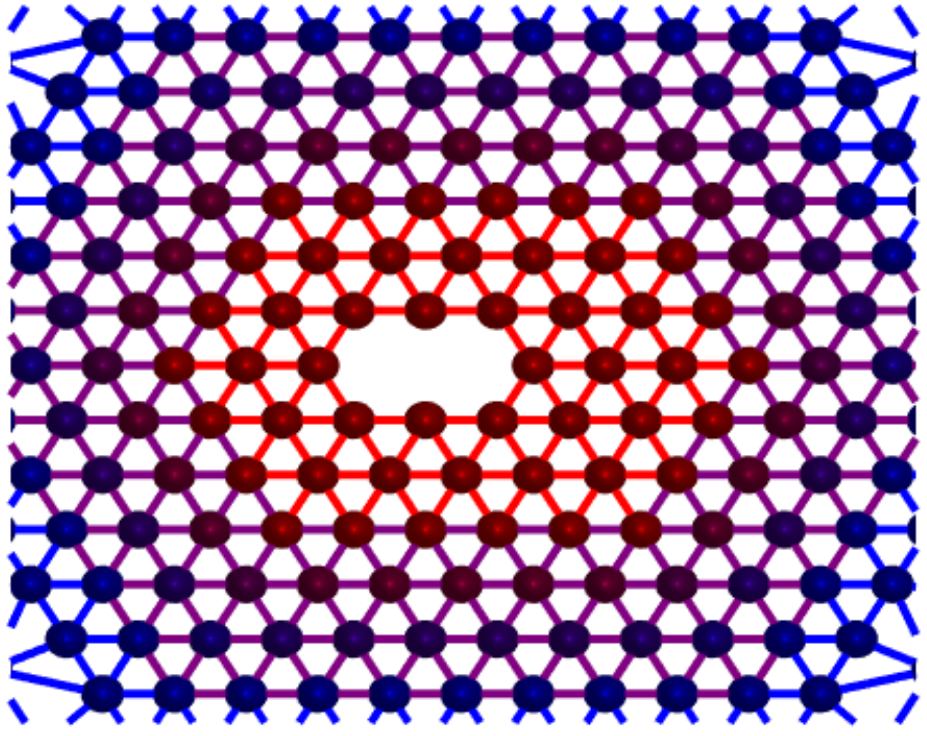}
	\caption{Computational domain, finite element grid and atomistic
region as used in the construction of blended a/c coupling methods. The color of the spheres in the right sub-figure represent the value of the
blending function (red as $\beta=0$ while blue as $\beta=1$).}
	\label{figs:geom_divac}
\end{center}
\end{figure}

Next, the computational domain $\Omega = \Omega^\a \cup \Omega^\b \cup \Omega^\c \subset \R^2$ can be decomposed into three regions, the {{\it atomistic region}} $\Omega^{\a}$ with radius $R^\a$, the {{\it blending region}} $\Omega^\b$ with width $L^{\b}$ and the {{\it continuum region}} $\Omega^\c$. We define the set of core atoms $\L^{\a} := \L \cap \Omega^\a$ and the set of blended atoms $\L^\b := \L \cap \Omega^\b$. Let $\T^{\rm ab}$ be the {\it canonical} background mesh induced by $\L^{\rm ab}:=\L^{\a}\cup\L^\b$, and $\T^{\rm c}_h$ be a shape-regular tetrahedral partition of the continuum region. We denote $\T_h = \T^{\rm ab} \bigcup \T^\c_h$ as the tetrahedral partition. The subscript $h$ represents the {\it coarse-graining} in the continuum region by applying the finite element methods. See Figure~\ref{figs:geom_divac} for an illustration of $\T_h$. The nodes of $\T_h$ is denoted as $\mathcal{X}_h$. Let ${\rm DoF}:=\#\mathcal{X}_h$ and $\Omega_{h} = \bigcup_{T\in \T_{h}} T$.

We consider a blending function $\beta \in C^{2,1}(\R^2)$ satisfying $\beta=0$ in $\Omega^{\a}$, $\beta=1$ in $\Omega^{\c}$ and ${\rm supp}(\nabla \beta) \subset \Omega^{\b}$. The choice of $\beta$ heavily influences the error estimates of the blended a/c coupling methods, we refer to \cite{2014-bqce, LuOr:acta, 2012-CMAME-optbqce} for a detailed discussion. In this work, following the construction in~\cite{2012-CMAME-optbqce, colz2016}, $\beta$ is obtained in a preprocessing step by approximately minimizing $\|\nabla^2 \beta\|_{L^2}$. 

If we use $\mathcal{P}_1$ finite element to discretize the Cauchy-Born continuum model, the space of {\it coarse-grained} displacements is then given by
\begin{eqnarray}
  \Us_{h}^{(1)} := \big\{ u_{h} \in C(\R^d;\R^d) ~\big|~
  \text{$u_{h} \in \mathcal{P}_1(\T_h)$, $u_{h} = 0$ in $\R^d \setminus \Omega_{h}$}\big\}.
\end{eqnarray}

To achieve the optimal rates of convergence in terms of the number of degrees of freedom ($\textrm{DoF}$) among all a/c coupling methods employing Cauchy-Born model in the continuum region~\cite[Theorem 2.1]{fang2020blended}, we can apply $\mathcal{P}_2$ finite element method in the continuum region.
Hence, we decompose $\T_h = \T_h^{\rm (1)}\cup\T_h^{\rm (2)}$, where
\[
\T_h^{\rm (1)} = \big\{T \in \T_h ~\big|~\beta|_{T}<1\big\}, 
\]
and replace $\Us_h$ with the approximate space
\begin{align}
    \Us_h^{(2)} := \big\{u_{h} \in C(\R^d;\R^d) ~\big|~
  \text{$u_{h} \in \mathcal{P}_1(\T^{\rm (1)}_h)$, $u_{h} \in \mathcal{P}_2(\T^{\rm (2)}_h)$, $u_{h} = 0$ in $\R^d \setminus \Omega_{h}$ }\big\}.
\end{align} 
That is, we retain the $\mathcal{P}_1$ discretization in the fully refined atomistic and blending region, but employ $\mathcal{P}_2$ finite elements in the continuum region.

\subsubsection{Cauchy-Born approximation}
\label{sec:sub:sub:cb}

To overwhelm the number of degrees of freedom and meanwhile preserve considerable accuracy, we first define a continuum approximation by the Cauchy-Born rule which is a typical choice in the multiscale context~\cite{e2007cb, OrtnerTheil2012}. The Cauchy-Born energy density functional $W : \R^{3 \times 3} \to \R$ reads
\begin{displaymath}
  W(\mF) := \det (\mA^{-1}) \cdot V(\mF \mathcal{R}),
\end{displaymath}
where $V$ is the homogeneous site energy potential on $\Lhom$.

We then define the pure Cauchy-Born finite element functional~\cite{2014-bqce}, 
\begin{eqnarray}
\E^{\rm cb}_h(u_h) := \int_{\Omega_h} Q_h\big[W(\nabla u_h)\big] \dx,
\end{eqnarray}
where $Q_h$ is the quadrature operator. In particular, it is the $\mathcal{P}_0$ midpoint interpolation operator when $\mathcal{P}_1$ finite element method is used in the continuum region, i.e., $u_h\in\Us^{(1)}_h$. While $u_h \in \Us_h^{(2)}$, for stability analysis~\cite[Section 4.1]{colz2016}, the operator $Q_h$ must now be adjusted to provide a third-order quadrature scheme (e.g., the face midpoint trapezoidal rule) so that $\nabla u_h \otimes \nabla u_h \in \Us_h^{(2)}$ can be integrated exactly.

\subsubsection{BQCE method}
\label{sec:sub:sub:bqce}

We define the BQCE energy functional for as
\begin{equation}\label{eq:Ebqce}
  \E_h^{\rm bqce}(u_{h}) := \sum_{\ell \in \L \cap \Omega_{h}}\big(1-\beta(\ell)\big)V_{\ell}\big(D u_h(\ell)\big) + \int_{\Omega_{h}} Q_{h} \big[ \beta(x) W(\D u_h) \big] \dx.
\end{equation}
The BQCE method is to find
\begin{align}\label{eq:variational-problem-BQCE}
u^{\rm bqce}_{h} \in& \arg\min \big\{ \E^{\rm bqce}_{h}(u_{h}),~u_{h} \in \Us_{h}^{(1)} \big\}. 
\end{align}
According to previous studies~\cite{colz2016, 2016-qcp2}, using higher-order finite elements in the blending region does not result in an improved BQCE method as the blending error dominates. Therefore, we do not take into account the $\mathcal{P}_2$ finite element when considering the BQCE method.

\subsubsection{BQCF method}
\label{sec:sub:sub:bqcf}

\def\F{\mathcal{F}}

While the BQCE method \eqref{eq:Ebqce} blends atomistic and continuum energies, the BQCF method \cite{2014-bqce, li2014theory} blends atomistic and continuum forces. It is worthwhile mentioning that there is no energy functional for this force-based method.

Recall that $\beta \in C^{2,1}(\R^2)$ is a blending function, and if we employ $\mathcal{P}_1$ finite element method, then the BQCF operator is the nonlinear map $\F^{\rm bqcf}_h:\Us^{(1)}_h \rightarrow (\Us^{(1)}_h)^{*}$, defined by
\begin{eqnarray}
\<\F^{\rm bqcf}_h(u_h), v_h\>:= \<\delta \E^{\rm a}(u_h), (1-\beta)v_h\> + \<\delta\E^{\rm cb}_h(u_h), \beta v_h\>. 
\end{eqnarray}
In the P1-BQCF method we solve the following variational nonlinear system
\begin{eqnarray}\label{eq:variational-problem-BQCF}
\<\F_h^{\rm bqcf}(u^{\rm bqcf,1}_{h}), v_h\> = 0, \quad \forall v_h \in \Us^{(1)}_h.
\end{eqnarray}
If one replaces $\Us_h^{(1)}$ by $\Us_h^{(2)}$, we obtain the P2-BQCF method as well as its solution $u^{\rm bqcf,2}_h$.

\subsubsection{BGFC method}
\label{sec:sub:sub:bgfc}

The energy functional of BGFC method is based on a {\it second renormalization} of the potential~\cite{colz2016} while the {\it first renormalization} is itself due to the assumption $V_{\ell}({\bf 0})=0$ made in Section~\ref{sec:sub:atom}. For $\ell \in \L$ and $u \in \Use$, we define
\begin{eqnarray}
V''_{\ell}\big(Du(\ell)\big):=V_{\ell}\big(Du(\ell)\big)-\<\delta V_{\ell}({\bf 0}), Du\>.
\end{eqnarray}
The corresponding second renormalized Cauchy-Born energy density is
\begin{eqnarray}
W''(\mF) := W(\mF) - \<\delta W({\bf 0}), \mF\>, \quad \text{for}~\mF\in \R^{d\times d}.
\end{eqnarray}
We then obtain the BGFC energy functional
\begin{eqnarray}\label{eqn:Eqc}
  \E^{\rm bgfc}_{h}(u_{h}) = \sum_{\ell \in \L \cap \Omega_{h}} (1-\beta(\ell))V''_\ell\big(Du_{h}(\ell)\big) + \int_{\Omega_{h}}Q_{h}\big[\beta(x) W''(\nabla u_{h}) \big]\dx.
\end{eqnarray}

Hence, the P1-BGFC and P2-BGFC methods are given respectively by
\begin{align}\label{eq:variational-problem-BGFC}
u^{\rm bgfc,1}_{h} \in \arg\min \big\{ \E^{\rm bgfc}_{h}(u_{h}), ~u_{h} \in \Us^{(1)}_{h} \big\}, \nonumber \\
u^{\rm bgfc,2}_{h} \in \arg\min \big\{ \E^{\rm bgfc}_{h}(u_{h}), ~u_{h} \in \Us^{(2)}_{h} \big\}.
\end{align} 

According to the discussions in \cite[Section 2.7]{colz2016}, the BGFC energy functional \eqref{eqn:Eqc} can be written as an equivalent ghost force removal formulation 
\begin{eqnarray}\label{eq:Ebgfc}
\E^{\rm bgfc}_{h}(u_{h}) = \E^{\rm bqce}_{h}(u_{h}) - \big\<\delta \E^{\rm bqce}_{\rm hom}({\bf 0}), u_{h}\big\>,
\end{eqnarray}
where $\E^{\rm bqce}_{\rm hom}({\bf 0})$ is the BQCE energy functional at the homogeneous lattice. We make fully use of this formulation throughout this paper due to its simplicity of implementation~\cite{fu2021adaptive, wang2022adapac2}.

\section{A Posteriori Error Estimates}
\label{sec:a_post}

\subsection{A priori error estimates}
\label{sec:sub:apriori}

The {\it a priori} error estimate for three typical blended a/c coupling methods (BQCE, BQCF, and BGFC methods) with both $\mathcal{P}_1$ and $\mathcal{P}_2$ finite element methods. The result is presented in Theorem~\ref{thm:ap} and is essentially derived from previous works such as \cite{colz2016, fang2020blended}. However, we adapt the result to the specific setting of this paper.

\begin{theorem}\label{thm:ap}
If $u^{\a}$ is a strongly stable solution of \eqref{eq:variational-problem} satisfying Assumption \ref{as:LS}. Under additional assumptions on the blending function and the mesh \cite[Assumption 1]{fang2020blended}, then there exist solutions to BQCE \eqref{eq:variational-problem-BQCE}, P1-(P2-)BQCF \eqref{eq:variational-problem-BQCF} and P1-(P2-)BGFC \eqref{eq:variational-problem-BGFC} methods, such that 
    \begin{align}\label{eq:conv_pd}
    &\| u^{\a} -  u^{\rm bqce}_{h}\|_{\UsH}
     \lesssim ({\rm DoF})^{1/2-2/d}, \nonumber \\
     \| u^{\a} - u^{\rm bqcf,1}_{h}\|_{\UsH} &\lesssim  ({\rm DoF})^{-1/2-1/d}, \qquad 
     \| u^{\a} - u^{\rm bqcf,2}_{h}\|_{\UsH} \lesssim  ({\rm DoF})^{-1/2-2/d}, \\
    \|u^{\a} -  u^{\rm bgfc,1}_{h}\|_{\UsH} &\lesssim ({\rm DoF})^{-1/2-1/d}, \qquad
    \|u^{\a} -  u^{\rm bgfc,2}_{h}\|_{\UsH} \lesssim ({\rm DoF})^{-1/2-2/d}. \nonumber
    \end{align}
Moreover, we have the corresponding energy error estimates
\begin{align}\label{eq:conv_E_pd}
    \big| \E^{\rm a}(u^{\a}) -  \E^{\rm bqce}_h(u^{\rm bqce}_{h})\big|
     &\lesssim ({\rm DoF})^{1-4/d}, \nonumber \\
    \big| \E^{\rm a}(u^{\a}) -  \E^{\rm bgfc}_h(u^{\rm bgfc,1}_{h})\big| &\lesssim ({\rm DoF})^{-1-2/d}, \\
    \big| \E^{\rm a}(u^{\a}) -  \E^{\rm bgfc}_h(u^{\rm bgfc,2}_{h})\big| &\lesssim ({\rm DoF})^{-1-2/d}.\nonumber
    \end{align}
\end{theorem}

\subsection{Residual-based a posteriori error estimate for geometry error}
\label{sec:Residual}

In the context of a/c coupling methods, a key challenge is to effectively determine the atomistic and continuum regions and to devise an appropriate mesh structure that strikes an optimal balance between accuracy and efficiency. To tackle this challenge, the {\it a posteriori} analysis and corresponding adaptive algorithms play a critical role in enabling the efficient implementation of a/c coupling methods.

In this section, we present the {\it a posteriori} error estimate for the geometry error. The estimate for the energy error will be examined in the subsequent section. Specifically, we consider the solution $u_h$ obtained from either the BQCE, P1-(P2-)BQCF, or P1-(P2-)BGFC model.
\[
u_h \in \{u^{\rm bqce}_h, u^{\rm bqcf,1}_h, u^{\rm bgfc,1}_h, u^{\rm bqcf,2}_h, u^{\rm bgfc,2}_h\}.
\]
Recalling that $I_{\a}u_h\in\Use$, we define the residual $\mR(I_{\a}u_h)$ as an operator on $\Use$, that is, 
\begin{equation}\label{eq:def_res}
\mR(I_{\rm a} u_{h})[v] := \< \del\Ea(\Ia u_{h}), v\>, \quad \forall v\in \Use.
\end{equation}

The following Lemma characterizes the dual norm of the residual and has been previously presented in our previous works~\cite{CMAME}. We include it here for the sake of completeness.

\begin{lemma}\label{lemma:res-F}
	Let $u^{\a}$ be a strongly stable solution of the atomistic model \eqref{eq:variational-problem}. If either BQCE, BQCF or BGFC is consistent in the sense of \eqref{eq:conv_pd}, then there exists a solution $u_{h}$ of \eqref{eq:variational-problem-BQCE} or \eqref{eq:variational-problem-BGFC}, and constants $c, C$ independent of the approximation parameters such that 
	\begin{eqnarray}\label{res-bound}
	c\|I_{\rm a} u_{h} - u^{\a}\|_{\UsH} \leq \| \mR(I_{\rm a} u_{h}) \|_{(\UsH)^*} \leq C\|I_{\rm a} u_{h} - u^{\a}\|_{\UsH},
	\end{eqnarray}
	where the residual $\textsf{R}(I_{\rm a} u_{h})$ is defined by \eqref{eq:def_res}.
\end{lemma}

Based on the previous Lemma, the optimal {\it a posteriori} error estimator for the geometry error can be derived as follows
\[
\eta^{\rm ideal}(u_{h}) := \| \mR(I_{\rm a} u_{h}) \|_{(\UsH)^*}.
\]
It is highly desirable as it offers both upper and lower bounds for geometry error up to independent constants. However, it cannot be computed in practice, which is a crucial limitation in the context of {\it a posteriori} error control and adaptive algorithms design.

\subsubsection{Residual stress-based error estimator}
\label{sec:sub:sub:imp_stress}

In this section, we present an alternative approach to residual estimates that exploits the ``implicit stress" of atomistic models, motivated by the {\it a posteriori} error estimates of QM/MM coupling methods \cite{CMAME, wang2020posteriori}.

Specifically, we consider the equivalent formulation of $\| \mR(I_{\a} u_{h}) \|_{(\UsH)^*}$ in terms of the solution of a whole-space auxiliary Poisson problem. This allows us to express the residual in terms of the {\it residual force} for any $v\in\UsH$.
\begin{equation}
	\label{eq:forcenew}
	\mR(I_{\a} u_{h})[v] = \< \delta\E^{\a}(I_{\rm a}u_{h}), v \> = \sum_{\ell \in \L} \F^{\rm a}_{\ell}(I_{\rm a}u_{h}) \cdot v(\ell),
\end{equation}
where the {\it residual force} $\F^{\rm a}_{\ell}(I_{\rm a}u_{h}) := \partial_{u(\ell)} \E^{\a}(u)\big|_{u = {I_{\rm a}u_{h}}}$. With this definition of residual force, we have the following Lemma and the proof is given in the Appendix~\ref{sec:app:analysis}.

\begin{lemma} \label{th:ctsphi}
	Up to a constant shift, there exists a unique $\phi(I_{\rm a} u_{h}; x) \in \dot{H}^1(\R^2)$ such that
	\begin{equation} \label{eq:defnphi}
	    \int_{\R^2} \nabla\phi(I_{\rm a} u_{h}; x) \cdot \nabla v {~\rm d}x
	    = \int_{\R^2} \widehat{\F}^{\rm a}(I_{\rm a} u_{h}; x) \cdot v {~\rm d}x \qquad \forall v \in H^1(\R^2),
	\end{equation}
	where 
	\begin{equation}
	\label{eq:rescaleIntp}
	\widehat{\F}^{\rm a}(I_{\rm a} u_{h}; x):= \sum_{\ell\in\L} c_{\ell} \F^{\rm a}_{\ell}(I_{\rm a} u_{h})\zeta_\ell(x) \qquad \textrm{with} \quad c_{\ell} = \frac{1}{\int_{\R^2} \zeta_\ell(x) {~\rm d}x}
	\end{equation}
	is the rescaled piecewise affine nodal interpolant of the residual forces.
	Moreover, there exist constants $c_1, C_1$ such that
	\begin{equation}
		c_1 \| \mR(I_{\rm a} u_{h}) \|_{(\UsH)^*} \leq \| \nabla \phi(I_{\rm a}u_{h}) \|_{L^2(\R^2)}
		\leq C_1 \| \mR(I_{\rm a} u_{h}) \|_{(\UsH)^*}.
		\label{eq:equiv}
	\end{equation}
\end{lemma} 

The presented Lemma yields several key observations: (1) the solution $\phi$ acts as a Riesz representation of the residual $\mR(I_{\a}u_{h})$; (2) the resulting equivalence demonstrates that $|\nabla\phi(I_{\rm a}u_{h})|_{L^2(\R^2)}$ can serve as an ideal error estimator since it can provide upper and lower bounds to the true approximation error; (3) $\nabla\phi$ corresponds to the second Piola atomistic stress tensor up to a divergence-free term~\cite{wang2018analysis}, which will be discussed further in Appendix~\ref{sec:app:Astress}. For these reasons, we refer to this error estimator as the residual stress-based error estimator.

\def\RO{R_{\Omega}}
\def\O{\Omega}

\subsubsection{Practical error estimator}

The estimator $|\nabla\phi|{L^2}$ obtained from Lemma~\ref{th:ctsphi} is considered an ideal error estimator because of the equivalence \eqref{eq:equiv}. However, it is not computable due to the implicit nature of equation \eqref{eq:defnphi}. To overcome this limitation, we propose a practical approach by (1) restricting the computational domain to a finite region $\Omega$ and (2) discretizing the Poisson problem using the canonical triangulation $\T_{\Omega}$ and the finite element method. 

Different from the discretization of Cauchy-Born approximation where $\mathcal{P}_2$ can be applied, it is sufficient to use $\mathcal{P}_1$ to solve \eqref{eq:defnphi} due to the discrete nature of atomistic stress. For simplicity, suppose that $\Omega$ is compatible with $\T_{\L}$, that is, there exists a subset $\T_{\Omega} \subset \T_{\L}$ such that ${\rm clos}(\Omega)=\cup\T_{\Omega}$. We then define the finite element space utilizing $\mathcal{P}_1$ finite element 
\[
	\mathscr{U}_{\a}^{\Omega}:=\{ \phi_{\a} \in \mathcal{P}_1(\T_{\Omega}): \phi_{\a} = 0\text{ on }\partial \Omega \}.
\]
We can now obtain an approximation to the idealised estimator $\phi$ by solving for $\phi_{\rm a} \in \mathscr{U}^{\Omega}_{\rm a}$ such that 
\begin{eqnarray}
\label{eq:stress_force_trun_T}
 \int_{\Omega}\nabla \phi_{\rm a}(I_{\a} u_{h}) \cdot \nabla v_{\a} {\dx} = \int_{\Omega} \widehat{\F}^{\a}(I_{\a} u_{h}) \cdot v_{\a} {\dx} \qquad \forall v_{\rm a} \in \mathscr{U}^{\Omega}_{\rm a}.
\end{eqnarray}
Hence, we define the approximate {\it a posterior} error estimator by
\begin{equation}\label{eq:APET}
	\eta(u_{h}) := \|\nabla \phi_{\a}(I_{\a} u_{h})\|_{L^2(\Omega)},
\end{equation}
and the indicator arising from the truncation by
\begin{align}\label{eq:rhoh_alt}
\rho_{\rm tr}(u_{h}) := R_{\Omega} \cdot \|\widehat{\mathcal{F}}^{\a}(I_{\a} u_{h})\|_{L^2(\Omega\setminus B_{R_\Omega/2})}
        +  \|\nabla \phi_\a(I_{\a} u_{h})\|_{L^2(\Omega\setminus B_{R_\Omega/2})}.
\end{align}
From the results in \cite{wang2020posteriori}, we can obtain the following theorem, demonstrating the equivalence of the {\it idealised residual estimate} $\|\mR(I_{\a} u_{h})\|_{(\Use)^*}$ and the {\it approximate estimator} $\eta(u_{h})$ up to the truncated indicator defined by \eqref{eq:rhoh_alt}. The proof is given in the Appendix~\ref{sec:app:analysis}.

\begin{theorem}\label{th:mainresult}
Let the approximate a posteriori error estimator $\eta(u_h)$ and 
the truncation indicator $\rho_{\rm tr}(u_h)$ be define by \eqref{eq:APET} and \eqref{eq:rhoh_alt} respectively. Under the conditions of Lemma \ref{lemma:res-F},	there exists positive constants $c, C$ such that
	\begin{align} \label{eq:main:first}
		c \eta(u_{h})
		\leq
		\| \mR(I_{\rm a} u_{h}) \|_{(\UsH)^*}
		\leq
		C \big( \eta(u_{h})
					+ \rho_{\rm tr}(u_{h}) \big).
	\end{align}
\end{theorem}

It is worthwhile to mention that compared to our previous work~\cite{wang2022adapac1}, the crucial improvements of the residual stress-based estimator are two-folds: (1) it provides both upper and lower bounds to the approximation error; (2) the resulting error from truncation is also specified rigorously. These will benefit to the design of corresponding adaptive algorithm given in Section~\ref{sec:adap}.  

\subsection{A posteriori error estimate for energy error}
\label{sec:sub:energy}

The total energy is a significant material property, in addition to the geometry error, in various applications. In this section, we aim to establish a rigorous {\it a posteriori} error estimate for the energy error, represented by $\big|\E^{\rm a}(u^{\a}) - \E^{\beta}{h}(u{h})\big|$, where $\E^{\beta}_{h}$ denotes the energy functional of one of the blended type a/c coupling methods. The analysis we present is fundamentally based on \cite{wang2018posteriori}, but we have modified it to suit the context of this work.

The energy error can be split into two parts
\[
\big|\E^{\rm a}(u^{\a}) - \E^{\beta}_{h}(u_{h})\big| \leq \big|\E^{\rm a}(u^{\a}) - \E^{\rm a}(I_{\a} u_{h})\big| + \big|\E^{\rm a}(I_{\a} u_{h}) - \E^{\beta}_{h}(u_{h})\big| =: E_1 + \big|\mu_{\rm E}(u_{h})\big|
\]
We first estimate the term $E_1$. As $\E^{\a}$ is twice differentiable along the segment $\{(1-s)u^{\rm a}+s I_{\a}u_{h}~|~ s\in(0,1)\}$ with a Lipschitz constant $M$, one can deduce from \eqref{eq:main:first} that 
\begin{align}
    |E_1| &= \Big| \int_{0}^{1} \<\delta \E^{\rm a}\big( (1-s)u^{\a}+s I_{\a}u_{h} \big), u^{\a}-u_{ h}\> \ds \Big| \nonumber \\
    &\leq M \|I_{\rm a} u_{h} - u^{\a}\|^2_{\Use} \leq C \big(\eta(u_{h})^2 + \rho_{\rm tr}(u_{h})^2\big).
\end{align}
For the second term $\mu_{\rm E}(u_{h})$, we utilize an equivalent version of atomistic energy functional such that it can be easily assigned into local contributions for both BQCE and BGFC methods. To that end, let $u_h^{\a}:=I_{\a}u_h$ for simplicity of notations and recall the definition \eqref{energy-difference}, we first rewrite $\E^{\rm a}$ in an element-wise style
\[
\E^{\rm a}(u^{\a}_{h}) = \sum_{T \in \T_{\L}} \frac{1}{6}\sum_{\ell \in T \cap \L} V_{\ell}\big(Du^{\a}_{h}(\ell)\big).
\]
For BQCE method, the energy functional can also be written as 
\begin{align*}
\E^{\rm bqce}_h(u_{h}) =& \sum_{T \in \T_{\L}} \frac{1}{6}\sum_{\ell \in T \cap \L} \big(1-\beta(\ell)\big)\cdot V_{\ell}\big(Du_{h}(\ell)\big) \\
&+\sum_{T\in\T_{h}}\sum_{T'\in\T_{\L}, T'\cap T\neq \emptyset} \beta_{T'} |T \cap T'| \cdot W(\nabla_{T'} u_{h}).
\end{align*}
Hence, $\mu_{\rm E}(u_h)$ for BQCE method reads 
\begin{align}\label{eq:muE_bqce}
    \mu_{\rm E}(u_{h}) =&~\E^{\rm a}(u^{\a}_{h}) - \E^{\beta}_{h}(u_{h}) \nonumber \\
=&  \sum_{T\in\T_{\rm b}} \Big( \frac{1}{6}\sum_{\ell\in T\cap \L} \beta(\ell) \cdot V_{\ell}\big(D u^{\a}_{h}(\ell)\big) - \beta_{T}|T|\cdot W(\nabla_{T}u_{h})\Big) \nonumber \\
    &+ \sum_{T \in \T_{\c}} \sum_{T' \in \T_{\L}, \omega(T')\cap \partial T \neq \emptyset} \frac{|T \cap T'|}{2|T'|}\Big( \frac{1}{3} \sum_{\ell \in T'\cap\L} V\big(Du^{\a}_{h}(\ell)\big) - V(\nabla u_{h}\cdot\Rg)\Big). 
\end{align}
The energy error estimate for BGFC method can be extended to include the correction term by replacing the potential functions $V$ and $W$ in \eqref{eq:muE_bqce} with $V''$ and $W''$, respectively, as per the definition in Section~\ref{sec:sub:sub:bgfc}.

\begin{theorem}\label{th:mainresultE}
Under the conditions of Lemma \ref{lemma:res-F}, let the approximate a posteriori error estimator $\eta(u_h)$ and 
the truncation indicator $\rho_{\rm tr}(u_h)$ be define by \eqref{eq:APET} and \eqref{eq:rhoh_alt} respectively, there exists positive constants $C_{\rm E}$ such that
\begin{align}
\big|\E^{\rm a}(u^{\a}) - \E^{\beta}_{h}(u_{h})\big| \leq~ C_{\rm E}\big(\eta(u_{h})^2 + \rho_{\rm tr}(u_{h})^2\big) + |\mu_{\rm E}(u_{h})|.
\end{align}
\end{theorem}
As a result, we define the energy estimator as 
\begin{eqnarray}\label{eq:etaE}
\eta_{\rm E}(u_h) := \eta(u_{h})^2 + |\mu_{\rm E}(u_{h})|.
\end{eqnarray}

\section{Adaptive Algorithms and Numerical Experiments}
\label{sec: numeric}

In this section, we first propose the residual-based local error estimators based on the theoretical results derived in the last section and design an adaptive algorithm for the {\it a posteriori} error control in Section~\ref{sec:adap} The adaptive computations for several typical examples of crystalline defects are then conducted in Section~\ref{sec: numexp}. 

\subsection{Adaptive algorithm}
\label{sec:adap} 

In order to devise and implement the adaptive algorithm, it is necessary to allocate the global residual estimator $\eta(y_h)$ into elementwise local contributions. Before introducing this, the role of the truncation indicator $\rho_{\rm tr}(u_h)$ is discussed. As a rule, the truncation residual estimator is calculated directly in practice because it is comparatively small if the computational domain $\Omega$ is adequately large. During the main adaptive process, it is checked whether it exceeds the threshold. It is important to note that the truncation residual is not assigned to local contributions, but rather computed directly using \eqref{eq:rhoh_alt}, and its value is checked during the adaptive process. This is because it is much smaller than the total local estimator $\rho_{T}$ when $\Omega$ is sufficiently large. The resulting Algorithm~\ref{alg:resF} provides additional details.

\subsubsection{Local error contributions}
\label{sec:sub:sub:local}

To design the adaptive algorithm that controls the extension of atomistic region in addition to the local mesh refinement in continuum region, we derive the local error estimators based on $\eta(u_h)$ defined in \eqref{eq:APET} as follows. 
 
For $T \in \T_{h}$, we define
\begin{equation}
    \label{eq:eta_local}
    \eta_{T}(u_h) := \sum_{T'\in \T_{\Omega}, T'\cap T \neq \emptyset} \frac{|T' \cap T|}{|T'|} \cdot \frac{\|\nabla \phi_\a(\Ia u_{h}) \|^2_{L^2(T')}}{\eta(u_h)}.
\end{equation}
It is straightforward to see $\sum_{T \in\T_{h}}\eta_{T}(u_h) = \eta(u_h)$. 

For the energy estimator $\mu_{\rm E}(u_h)$ it is already given in~\eqref{eq:etaE}, we can define the local contributions similarly as $\mu_{\rm E}(T; u_h)$ such that $\sum_{T \in\T_{h}}\mu^2_{\rm E}(T; u_h) = \mu^2_{\rm E}(u_h)$. Meanwhile, for the energy based estimate, we have, 
\begin{equation}
    \label{eq:rhoexact_E}
    \eta^{\rm E}_{T}(u_h) := \sum_{T'\in \T_{\Omega}, T'\cap T \neq \emptyset} \frac{|T' \cap T|}{|T'|} \cdot \|\nabla \phi_\a(\Ia u_{h}) \|^2_{L^2(T')} + |\mu_{\rm E}(T; u_h)|.
\end{equation}
The additional energy part $\mu_{\rm E}(u_h)$ is already defined by elements, and therefore it is straightforward to get $\mu_{\rm E}(T;u_h)$. Notice that the sum of local estimators is equal to the global estimator.

\subsubsection{Adaptive algorithm}
\label{sec:sub:sub:alg}

Our adaptive procedure follows from the classic one~\cite{wang2018analysis, wang2020posteriori, wang2022adapac1}, namely
$$
\mbox{Solve}~\rightarrow~\mbox{Estimate}~\rightarrow~
\mbox{Mark}~\rightarrow~\mbox{Refine}.
$$
The main adaptive algorithm essentially follows from \cite[Algorithm 3]{wang2018posteriori} but we adapt it to the current setting of blended a/c coupling methods, which is summarized below. We are ready to give the main adaptive algorithm of the error estimators based on the implicit stresses. 

\begin{algorithm}
\caption{Adaptive blended a/c coupling based on residual forces}
\label{alg:resF}
\begin{enumerate}
	\item[Step 0] Prescible $\Omega$, $\T_h$, $N_{\max}$, $\eta_{\rm tol}$, $\tau_1$ and $\tau_2$.
	\item[Step 1] \textit{Solve}: Solve the BQCE (BQCF or BGFC) solution $u_{h}$ of \eqref{eq:variational-problem-BQCE} (\eqref{eq:variational-problem-BQCF} or \eqref{eq:variational-problem-BGFC}) on the current mesh $\T_{h}$.  
	\item[Step 2] \textit{Estimate}: Compute the local error estimator $\eta_{T}(u_h)$ by \eqref{eq:eta_local} for each $T\in\T_h$, and the truncation indicator $\rho_{\rm tr}(u_h)$ by \eqref{eq:rhoh_alt}. Compute the degrees of freedom $N$ and $\tilde{\eta}^{\O}:=\sum_{T} \tilde{\eta}^{\O}_{T}$. If $\rho^{\rm tr} > \tau_1 \tilde{\eta}^{\O}$, enlarge the computational domain. Stop if $N>N_{\max}$ or $\tilde{\eta}^{\O} < \eta_{\rm tol}$. 
	\item[Step 3] \textit{Mark}: Apply Algorithm \ref{alg:mark} to construct the refinement set $\M$, the number of atomistic layers $k$ and the ratio $\alpha$.
	\item[Step 4] \textit{Refine:} Expand the atomistic region outward by $[\alpha k]$ layers and the blending region outward by $k-[\alpha k]$ layers. Construct $\M_{h}:=\{T\in \T_{h}: \exists T'\in \M, T \cap T' = T'\}$. Bisect all elements $T\in \M_{h}$. Go to Step 1.
\end{enumerate}
\end{algorithm}

We aim to provide a comparison between our main adaptive algorithm (Algorithm \ref{alg:resF}) and \cite[Algorithm 3]{wang2018posteriori}. We note that \cite[Algorithm 3]{wang2018posteriori} focuses on enlarging the atomistic region, while our algorithm extends both the atomistic and blending regions during the adaptive procedure.

Toward that end, we first introduce some notations. Let ${\rm dist}(T, \L_{\a}) := \inf\{|\ell-x_T|,\forall\ell\in\L_{\a}\}$ be the distance between $T$ and $\L_{\a}$ with $x_T$ the barycenter of $T$. Similarly we can define ${\rm dist}(T, \L_{\c})$ and ${\rm dist}(T, \L_{\a}\cup\L_{\b})$. The {{\it Mark}} step mainly applies the standard D\"{o}rfler strategy \cite{Dorfler1996} to choose a set $\M$ for model refinement and then we design a strategy demonstrating how to enlarge $\L_{\a}$ and $\L_{\b}$ based on the distance to different regions. We summarize it as follows, which is mainly based on the algorithm proposed in our previous work~\cite{wang2022adapac2}.

\begin{algorithm}
\caption{Mark step}
\label{alg:mark}
\begin{enumerate}
	\item[Step 1]:  Choose a minimal subset $\M\subset \T_{h}$ such that
	\begin{displaymath}
		\sum_{T\in\M}\tilde{\eta}^{\O}_{T}\geq\frac{1}{2} \tilde{\eta}^{\Omega}.
	\end{displaymath}	 
	\item[Step 2]: We can find the interface elements which are within $k$ layers of atomistic distance, $\M^k:=\{T\in\M\bigcap (\T_{\b} \cup \T_{\c}): {\rm{dist}}(T, \Lambda_{\a})\leq k \}$. Choose $K\geq 1$, find the first $k\leq K$ such that 
	\begin{equation}
		\sum_{T\in \M^k}\tilde{\eta}^{\O}_{T}\geq \tau_2\sum_{T\in\M}\tilde{\eta}^{\O}_{T},
	\end{equation}
	with tolerance $0<\tau_2<1$. 
	\item[Step 3]: Divide $\M^k:=\M^k_{\a} \cup \M^k_{\c}$ into two adjacent subsets such that $\M^k_{\a}:=\{T\in\M^k\cap\T_{\b}: {\rm dist}(T, \L_{\a}) \leq \theta L^{\b}{\rm dist}(T, \L_{\c})\}$. Compute the ratio 
	\begin{eqnarray}\label{eq:ratio}
	\alpha:=\frac{\sum_{T\in\M^k_{\a}} \tilde{\eta}^{\O}_T}{\sum_{T\in \M^k}\tilde{\eta}^{\O}_{T}}.
	\end{eqnarray}
	Let $\M:=\M\setminus\M^k$.
\end{enumerate}
\end{algorithm}

\begin{remark}
The main contribution of Algorithm~\ref{alg:mark} is that it introduces a control parameter $\theta$ inside the mesh refinement step to gain more information from the error estimator belonging to different regions. This allows for a more balanced refinement between the atomistic and blending regions, resulting in a more accurate and efficient adaptive algorithm. The parameter $\theta$ is defined as $\T^{\rm ref}{\b}(\theta):={T\in\T{\b}: {\rm dist}(T, \L_{\a}) \leq \theta L^{\b}}$, where $L^{\b}$ is the size of the blending region. By adjusting $\theta$, one can control the number of elements in the blending region to achieve a desired accuracy.

To determine the optimal value of $\theta$, we use a bisection method to find the value that gives an approximately equal total error contribution between $\T^{\rm ref}{\b}(\theta)$ and $\T{\b}$. However, we note that this approach can be further improved by considering anisotropic refinement of the interface, which can be achieved by evolving the interface instead of just updating the radius of each sub-region. We discuss this and other potential future developments in Section~\ref{sec: conclusion}.
\end{remark}

\subsection{Numerical experiments}
\label{sec: numexp}

In this section, based on our main adaptive algorithm proposed in Algorithm~\ref{alg:resF}, we conduct the adaptive computations for two typical types of crystalline defects, namely micro-crack and Frenkel defect.

\subsubsection{Model setup}
\label{sec:sub:sub:model}
Our prototype implementation for three blended a/c coupling methods is the 2D triangular lattice $\Lhom = \mathsf{A}\Z^2$ defined by
\begin{displaymath}
{\sf A}=\left(
	\begin{array}{cc}
		1 & \cos(\pi/3) \\
		0 & \sin(\pi/3)
	\end{array}	 \right),
\end{displaymath}
where $\Lhom$ is in fact the projection of a BCC crystal along the (111) direction.
%

We use the well-known EAM potential as the site potential $V_{\ell}$ through our numerical experiments, where
\begin{align}\label{eq:eam_potential}
  V_{\ell}(y) := & \sum_{\ell' \in \Nhd_{\ell}} \phi(|y(\ell)-y(\ell')|) + F\B(
  {\textstyle \sum_{\ell' \in \Nhd_{\ell}} \psi(|y(\ell)-y(\ell')|)} \B),\nonumber\\
    = &\sum_{\rho \in \Rg_{\ell}} \phi\big(|D_\rho y(\ell)|\big) + F\B(
  {\textstyle \sum_{\rho \in \Rg_{\ell}}} \psi\big( |D_\rho y(\ell)|\big) \B),  
\end{align}
for a pair potential $\phi$, an electron density function $\psi$ and
an embedding function $F$. In particular, we choose
\begin{displaymath}
\phi(r)=\exp(-2a(r-1))-2\exp(-a(r-1)),\quad \psi(r)=\exp(-br)
\end{displaymath}
\begin{displaymath}
F(\tilde{\rho})=C\left[(\tilde{\rho}-\tilde{\rho_{0}})^{2}+
(\tilde{\rho}-\tilde{\rho_{0}})^{4}\right]
\end{displaymath}
with parameters $a=4, b=3, c=10$ and $\tilde{\rho_{0}}=6\exp(0.9b)$, which is the same as the numerical experiments in \cite{wang2018posteriori, liao2018posteriori}. Compared to these previous works, one of the main contribution of this work is that we admit the finite range interactions. We typically choose $r_{\rm cut}=2$, that is, the next nearest neighbor interactions are taken into consideration for our numerical experiments. We also note that the specific choice of site potential is not essential and in fact we expect that most common site potentials would lead to similar results.

We mainly follow the constructions of blended a/c coupling methods in \cite{colz2016}, where the blending function is obtained in a preprocessing step by approximately minimising $\|\nabla^2 \beta\|_{L^2}$, as described in detail in \cite{2012-CMAME-optbqce}. For BGFC method, we implement the equivalent ``ghost force removal formulation" \eqref{eq:Ebgfc} instead of the ``renormalisation formulation" \eqref{eqn:Eqc} since it is relatively easy to be implemented.

We numerically test the main adaptive algorithm (Algorithm~\ref{alg:resF}) with the {\it approximated} error estimator defined by \eqref{eq:eta_local}. To better reveal the effectiveness and the efficiency of our error estimators, we compare the results of our adaptive computations with those using an {\it a priori} graded mesh \cite{2012-CMAME-optbqce, 2014-bqce}. In all numerical experiments, the radius of computational domain $\Omega$ is 300 initially. The adaptive processes start with the initial configuration with $R_{\rm a}=3$. The adaptive parameters in Algorithm~\ref{alg:resF} are fixed to be $\tau_1=1.0$ and $\tau_2=0.7$.

\subsubsection{Micro-crack}
\label{sec:sub:sub:mcrack}

We apply an isotropic stretch $\mathrm{S}$ and shear $\gamma_{II}$ by setting
\begin{displaymath}
{\sf B}=\left(
	\begin{array}{cc}
		1 & \gamma_{II} \\
		0            & 1+\mathrm{S}
	\end{array}	 \right)
	\cdot{\sf {F_{0}}}
\end{displaymath}
where ${\sf F_{0}} \propto \mathrm{I}$ minimizing the Cauchy-Born energy density $\mathrm{W}$ and $\mathrm{S}=\gamma_{II}=0.03$.

\begin{figure}[htb]
\begin{center}
	\includegraphics[height=6cm]{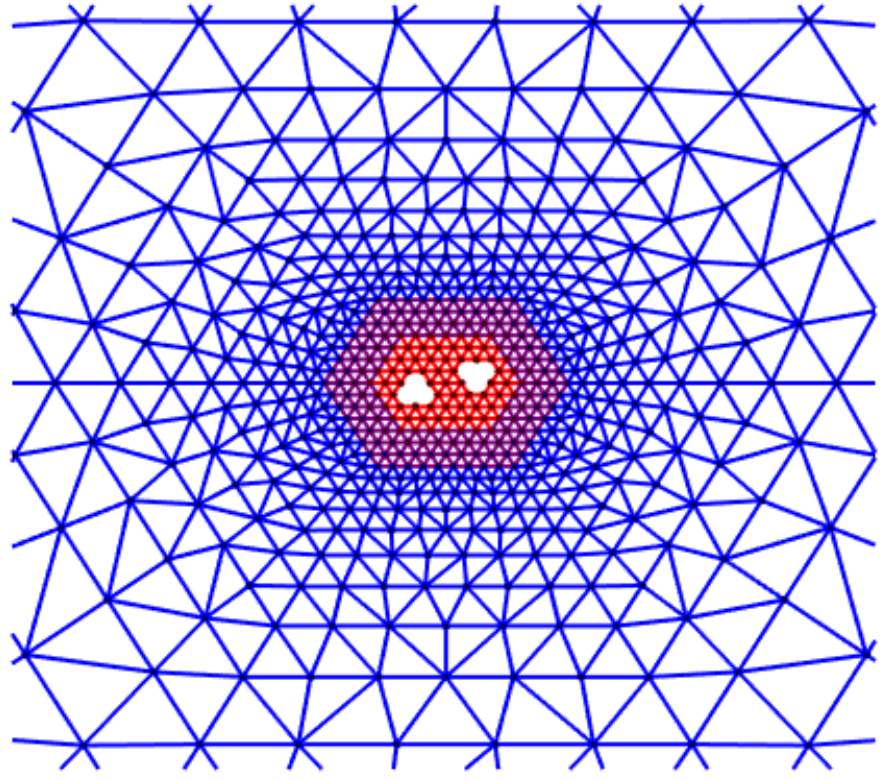}\quad
 \includegraphics[height=6cm]{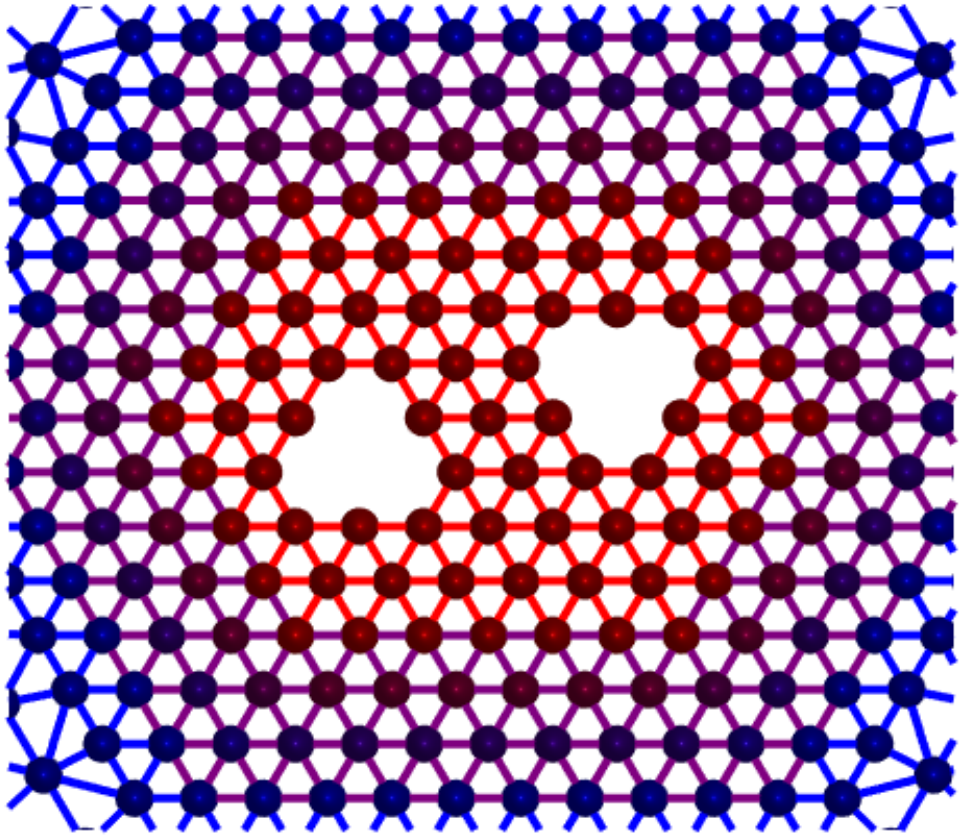}
	\caption{Computational domain, finite element grid and atomistic
region as used in the construction of blended a/c coupling methods for micro-crack.}
	\label{figs:geom_mcrack}
\end{center}
\end{figure}

The first defect we consider is the micro-crack, which is a prototypical example of point defects. We note that this is not technically a ``crack", it serves as an example of a localized defect with an anisotropic shape. To generate this defect, we remove $k$ atoms from $\Lhom$. In practice, we set $k=6$. The geometry of micro-crack is illustrated in Figure \ref{figs:geom_mcrack}.

\begin{figure}[htb]
\begin{center}
	\includegraphics[height=6cm]{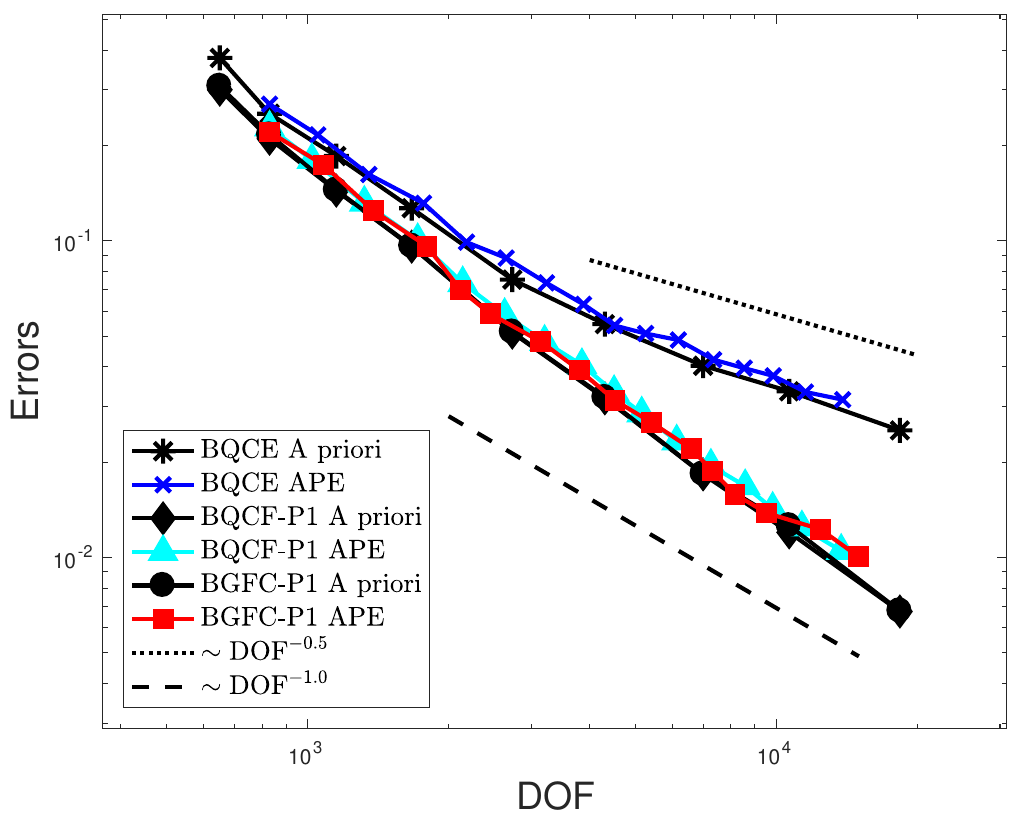}
 \includegraphics[height=6cm]{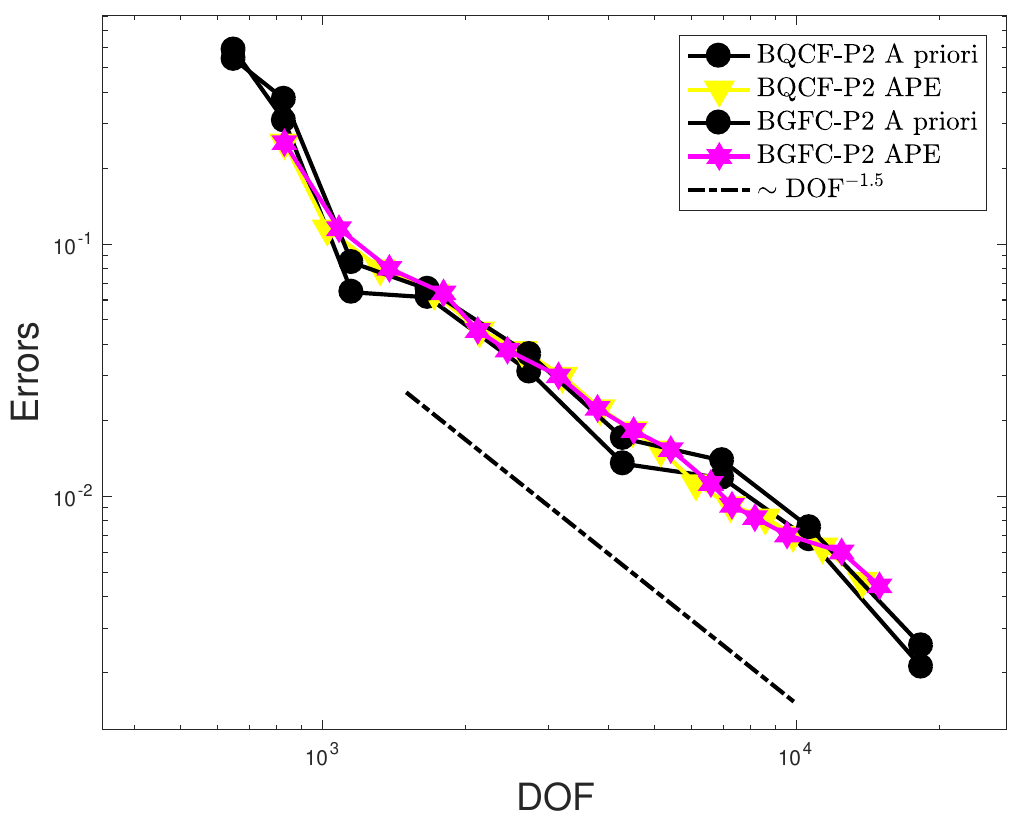}
	\caption{The convergence of geometry error for various a/c coupling methods with respect to the number of degrees of freedom for micro-crack.}
	\label{fig:conv_mcrack_conv}
\end{center}
\end{figure}

We observe in Figure~\ref{fig:conv_mcrack_conv} that all adaptive computations achieve optimal convergence rate ($N^{-0.5}$ for BQCE, $N^{-1.0}$ for BQCF-P1 and BGFC-P1, and $N^{-1.5}$ for BQCF-P2 and BGFC-P2) compared with the {\it a priori} graded mesh given by \cite{colz2016, fang2020blended}. The accuracy and robustness of the main adaptive algorithm are also demonstrated in this case.

\begin{figure}[htb]
\begin{center}
	\includegraphics[height=6cm]{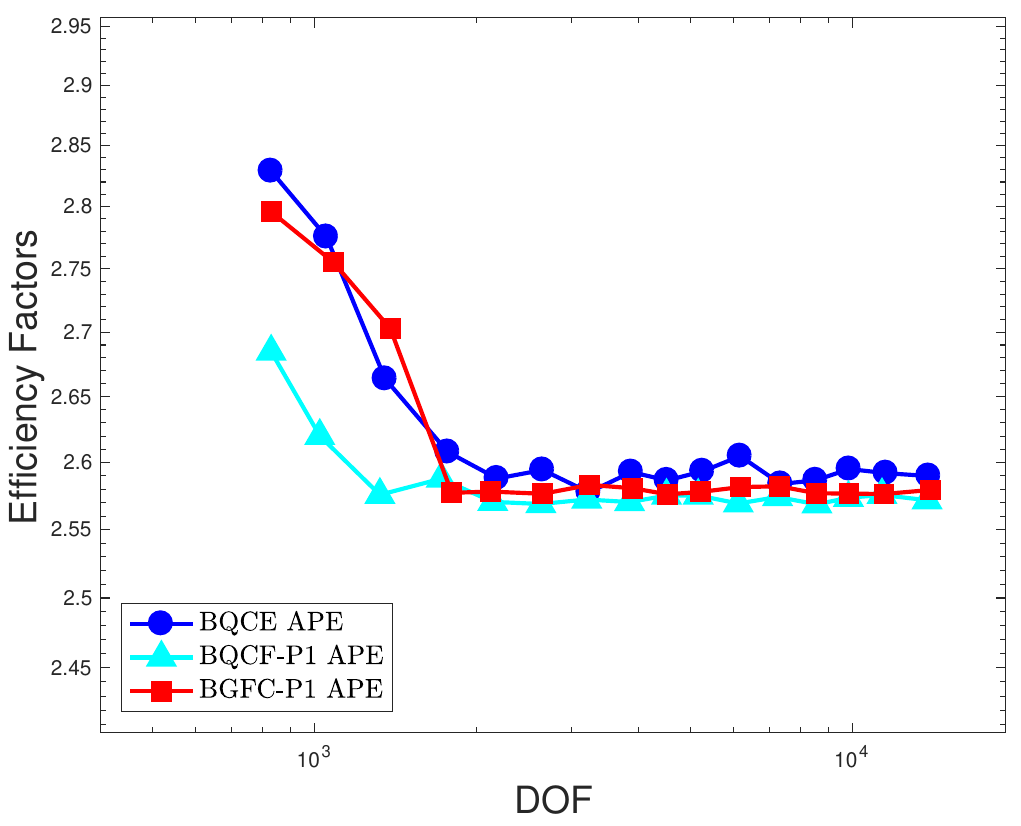}
        \includegraphics[height=6cm]{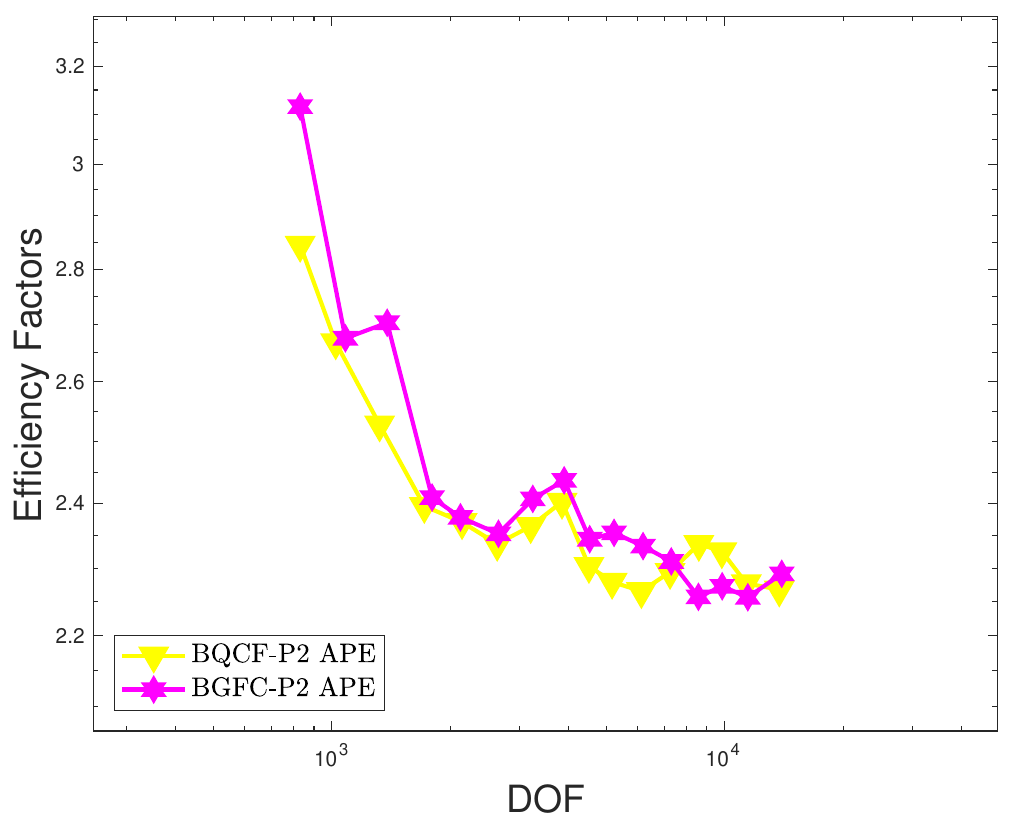}
	\caption{The efficiency factors of error estimators for different a/c coupling methods with respect to the number of degrees of freedom for micro-crack.}
	\label{fig:conv_mcrack_efffac}
\end{center}
\end{figure}

Figure~\ref{fig:conv_mcrack_efffac} plots the efficiency factors for both {\it rigorous} and {\it approximated} error estimators for the blended a/c coupling methods. We observe that the efficiency factors are all moderate during the adaptive computations for all a/c coupling methods, which is a novel improvement compared to the previous works for residual-force based error estimators~\cite{wang2022adapac2}. This provides a numerical evidence of proving both reliability and efficiency for the proposed residual-stress based error estimator.

\begin{figure}[htb]
\begin{center}
	\includegraphics[height=6cm]{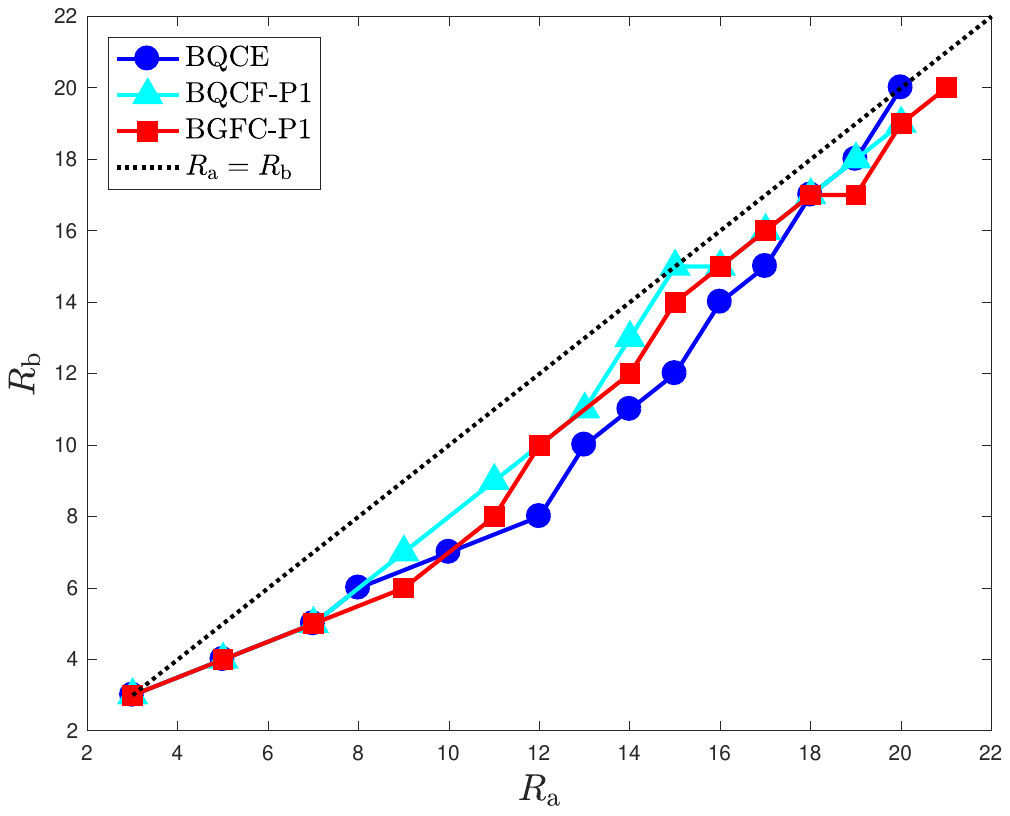}
 \includegraphics[height=6cm]{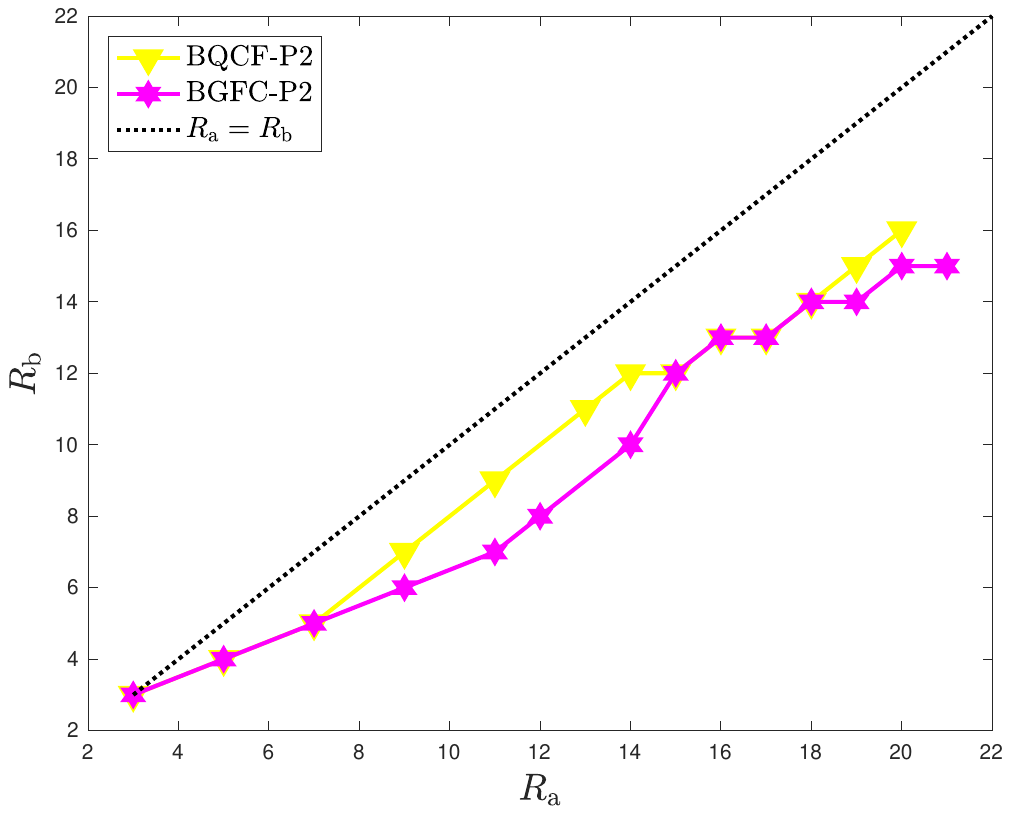}
	\caption{The relationship between the radius of the atomistic region $R_{\a}$ and the width of the blending region $R_{\b}$ in the adaptive computations for micro-crack.}
	\label{figs:RaRb_mcrack}
\end{center}
\end{figure}

We plot the relationship between $R_{\a}$ and $R_{\b}$ in the adaptive computations for micro-crack in Figure~\ref{figs:RaRb_mcrack}, where $R_{\a}$ is the radius of the atomistic region while $R_{\b}$ represent the width of the blending region. Our adaptive algorithm can achieve the optimal relationship between $R_{\a}$ and $R_{\b}$ for the case of micro-crack, indicating the robustness of main adaptive algorithm.

\begin{figure}[htb]
\begin{center}
	\includegraphics[height=6cm]{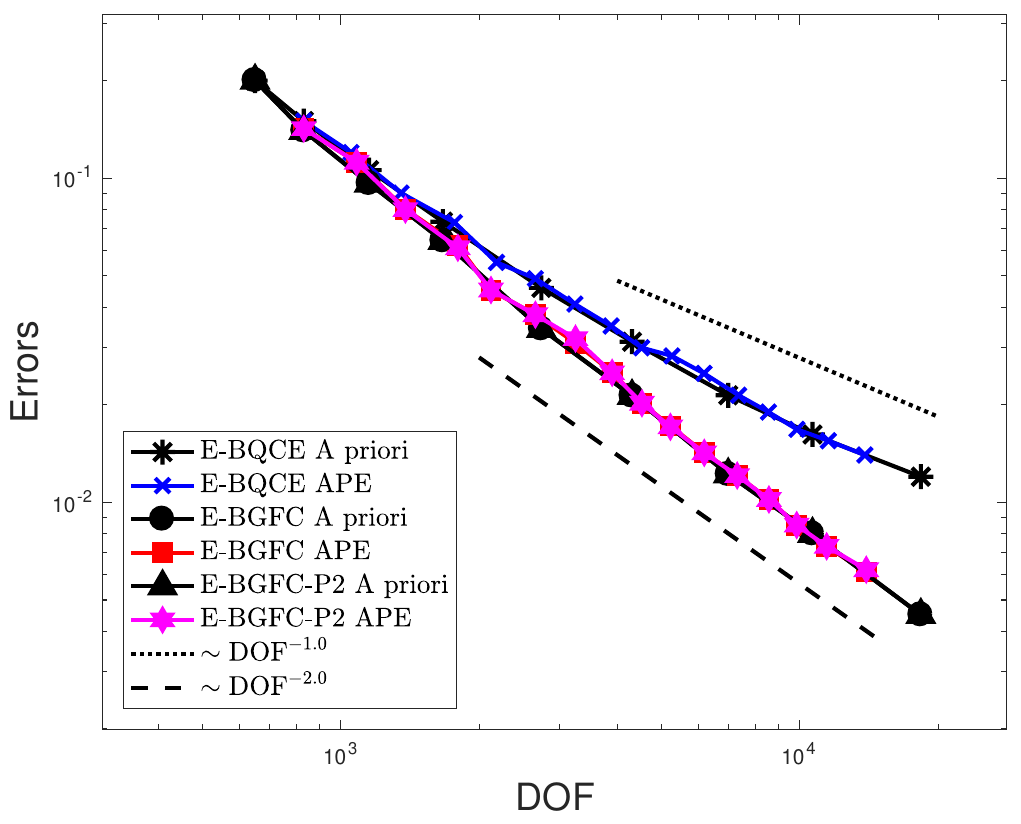}
	\caption{The convergences of energy error for both a priori and a posteriori errors with respect to the number of degrees of freedom for micro-crack.}
	\label{fig:conv_E_mcrack_conv}
\end{center}
\end{figure}

As for the energy error, we observe in Figure~\ref{fig:conv_E_mcrack_conv} that all adaptive computations achieve the same optimal convergence rate ($N^{-1.0}$ for BQCE whereas $N^{-2.0}$ for BGFC-P1 and BGFC-P2) compared with the {\it a priori} graded mesh given by \cite{colz2016, fang2020blended}, which indicate the error estimator for energy error is still reliable.

\subsubsection{Frenkel defect}
\label{sec:sub:sub:inters}

\begin{figure}[htb]
\begin{center}
	\includegraphics[height=6cm]{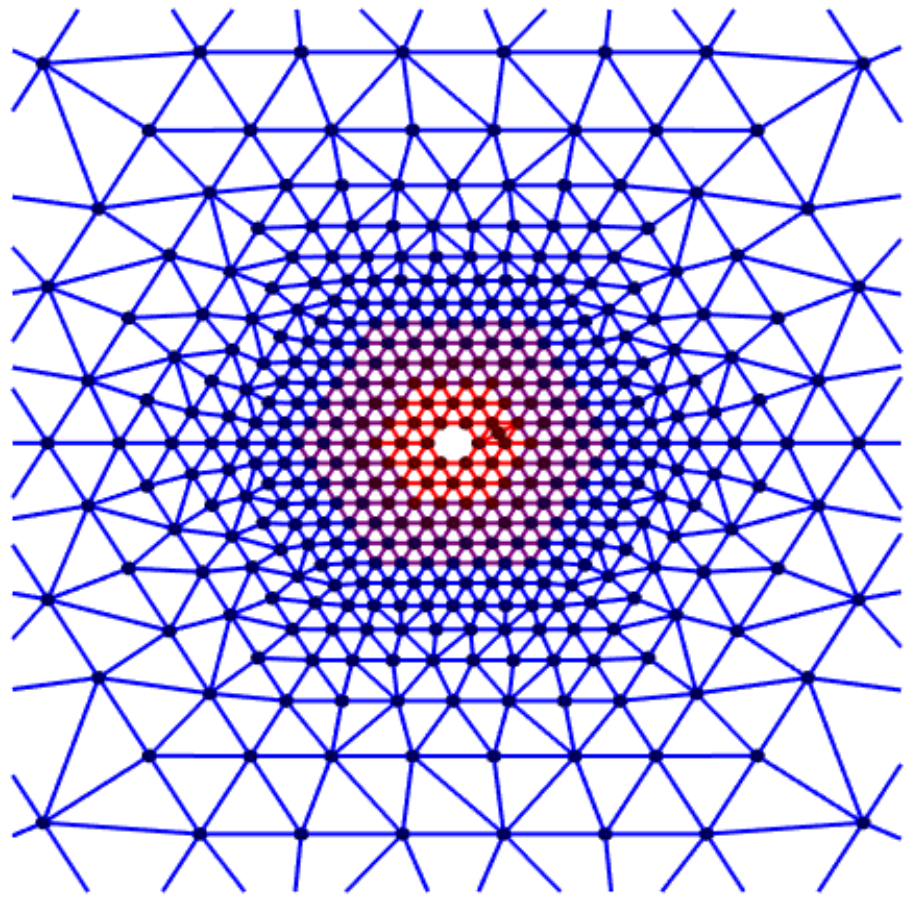}\quad
 \includegraphics[height=6cm]{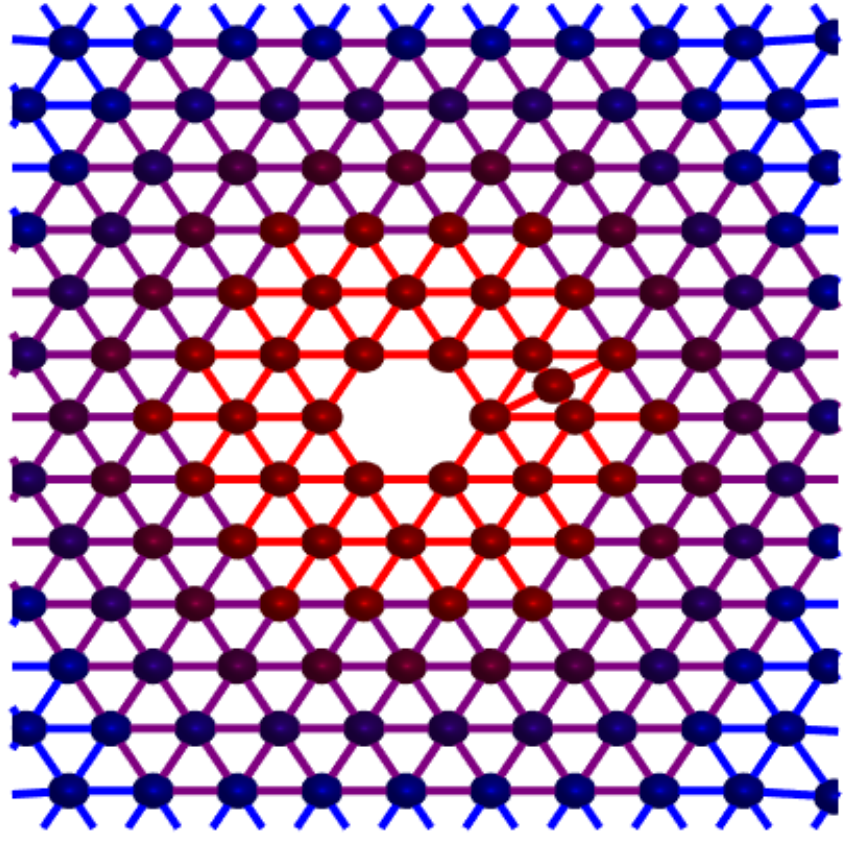}
	\caption{Computational domain, finite element grid and atomistic
region as used in the construction of blended a/c coupling methods for Frenkel defect.}
	\label{figs:geom_frenkel}
\end{center}
\end{figure}

In our second numerical experiment, we consider the well-known Frenkel defect to justify the performance of the main algorithm. We consider the lattice with one vacancy and one interstitial that are close, in this case, $\L= (\Lhom \setminus \L_1^{\rm def})\cup \{(3/2, 0)\}$. This is a very common point defect type but has never been considered so far in the literature of adaptive a/c methods to the best knowledge of the author. The geometry of micro-crack is illustrated in Figure \ref{figs:geom_frenkel}. We apply the same isotropic stretch $\mathrm{S}$ and shear $\gamma_{II}$ as that in previous example.

\begin{figure}[htb]
\begin{center}
	\includegraphics[height=6cm]{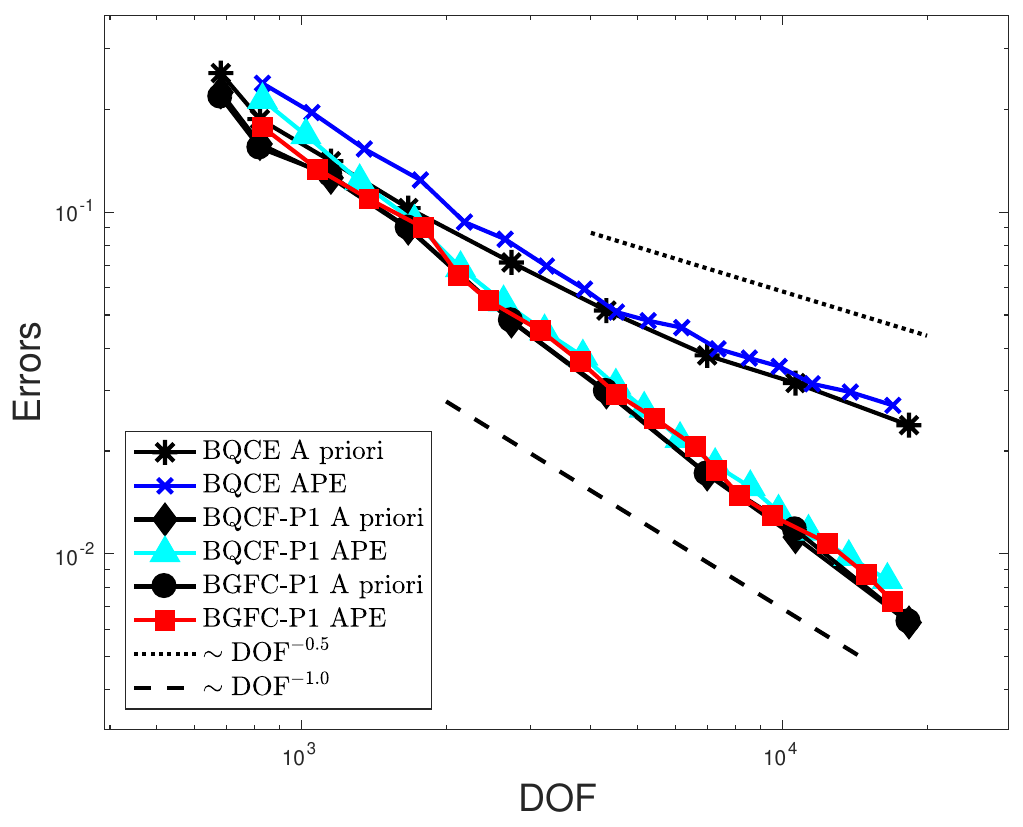}
 \includegraphics[height=6cm]{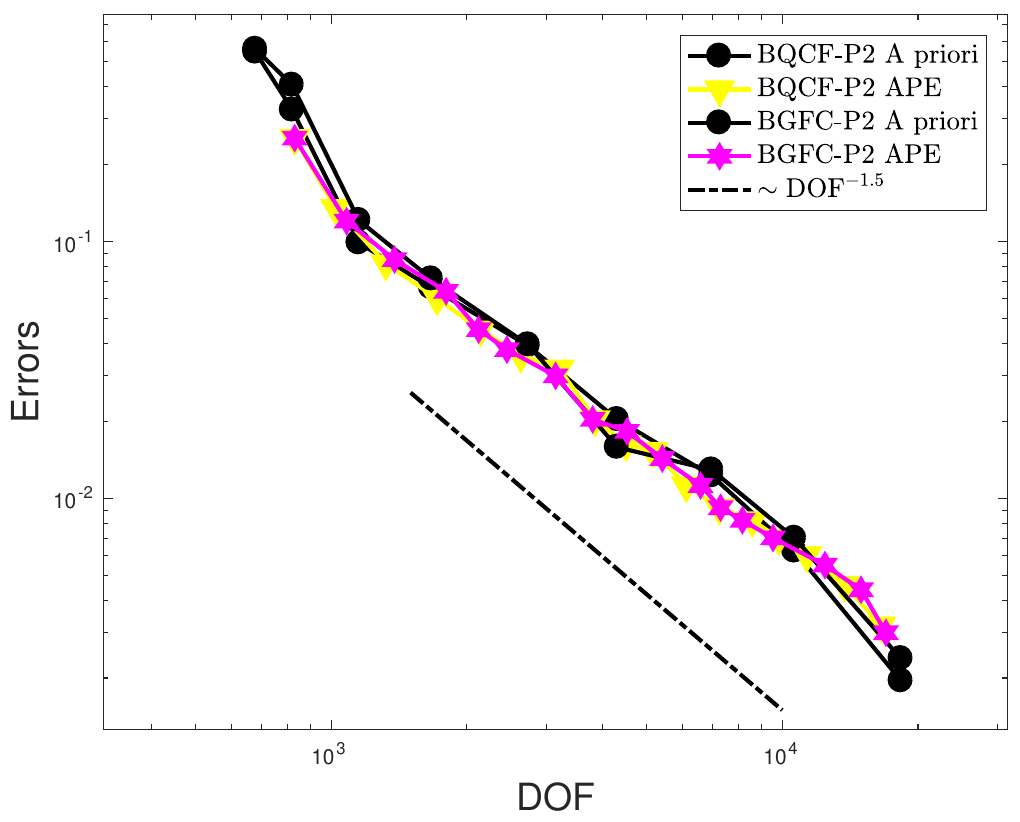}
	\caption{The convergences of geometry errors with respect to the number of degrees of freedom for Frenkel defect.}
	\label{fig:conv_frenkel_conv}
\end{center}
\end{figure}

We observe in Figure~\ref{fig:conv_frenkel_conv} that all adaptive computations achieve optimal convergence rate ($N^{-0.5}$ for BQCE, $N^{-1.0}$ for BQCF-P1 and BGFC-P1, and $N^{-1.5}$ for BQCF-P2 and BGFC-P2) compared with the {\it a priori} graded mesh given by \cite{colz2016, fang2020blended}. The accuracy and robustness of the main adaptive algorithm are also demonstrated in this case.

\begin{figure}[htb]
\begin{center}
	\includegraphics[height=6cm]{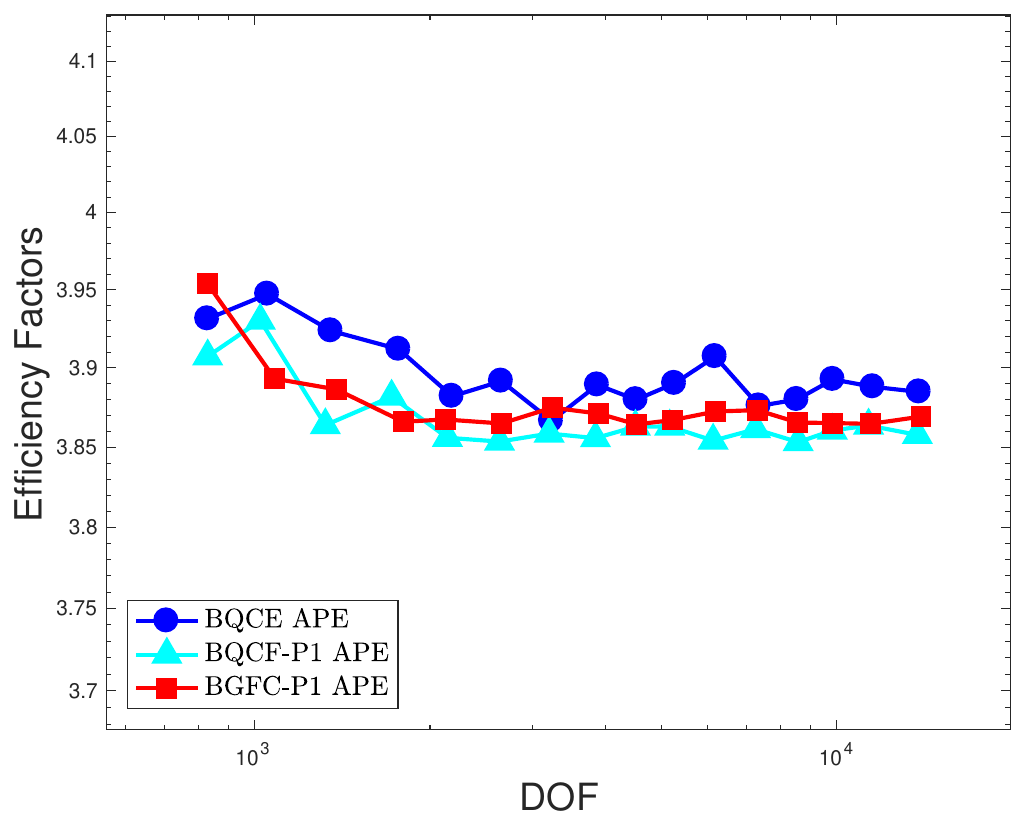}
 \includegraphics[height=6cm]{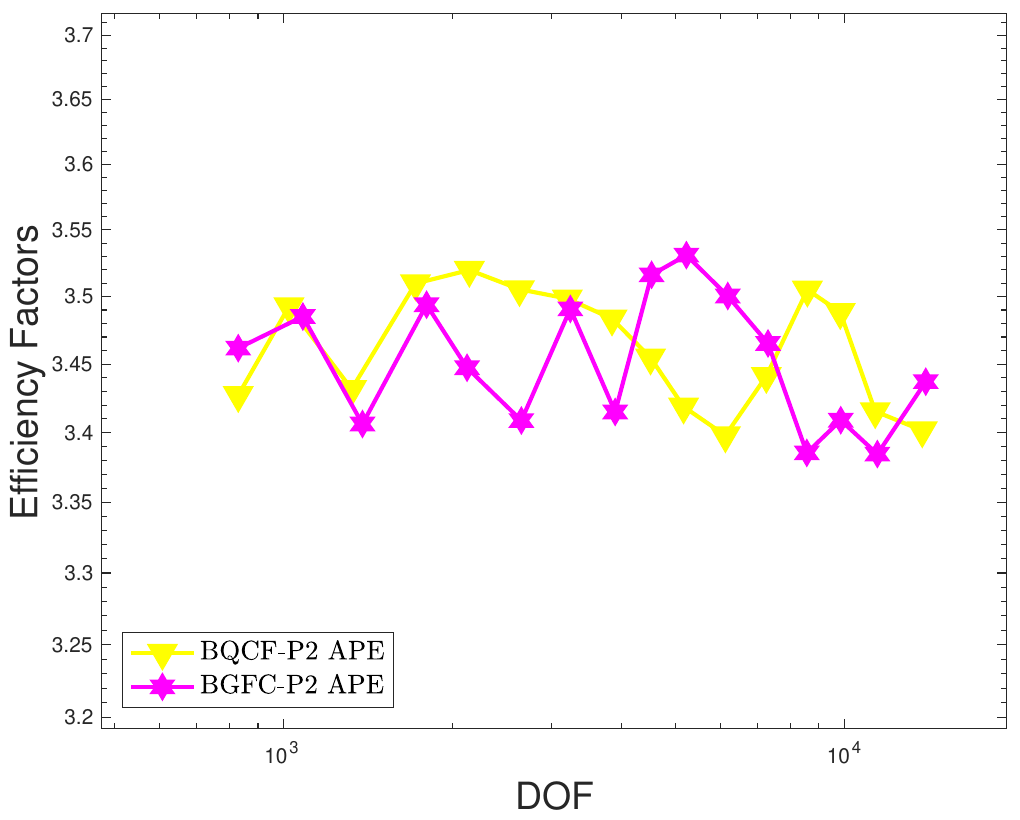}
	\caption{The efficiency factors of error estimator for various a/c coupling methods with respect to the number of degrees of freedom for Frenkel defect.}
	\label{fig:conv_frenkel_efffac}
\end{center}
\end{figure}

Figure~\ref{fig:conv_frenkel_efffac} plots the efficiency factors for both {\it rigorous} and {\it approximated} error estimators for the blended a/c coupling methods. We observe that the efficiency factors are all moderate during the adaptive computations for all a/c coupling methods, which is a novel improvement compared to the previous works for residual-force based error estimators~\cite{wang2022adapac2}. This provides a numerical evidence of proving both reliability and efficiency for the proposed residual-stress based error estimator.

\begin{figure}[htb]
\begin{center}
	\includegraphics[height=6cm]{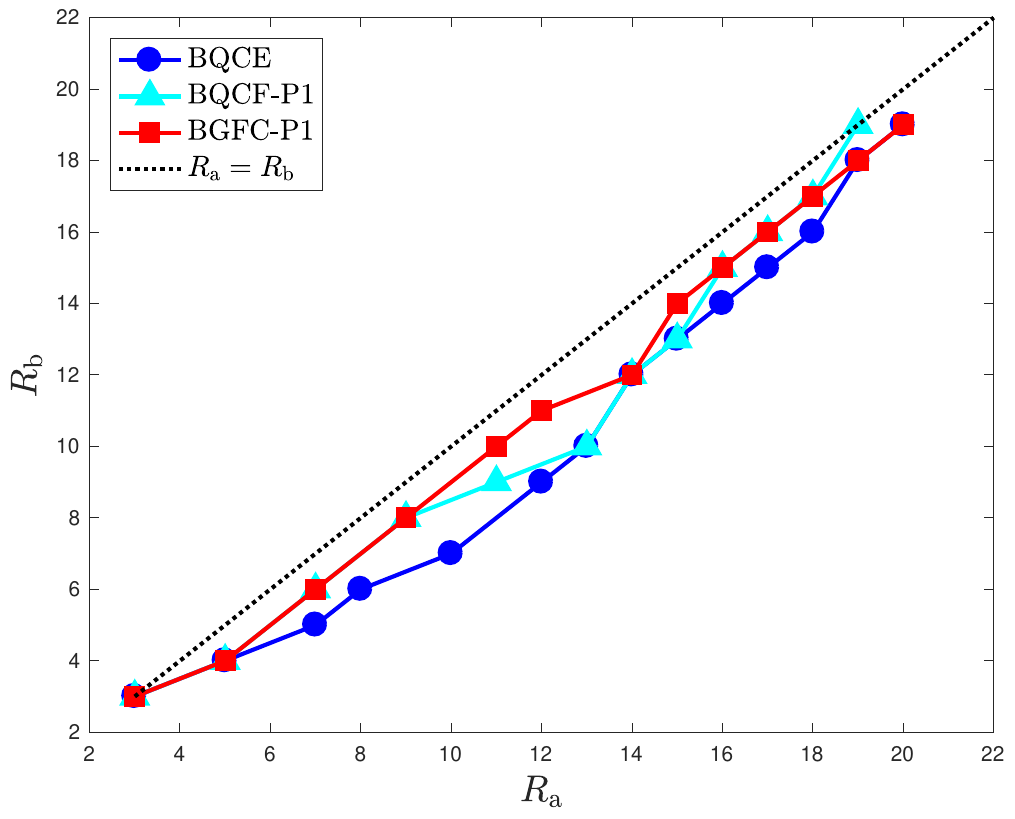}
 \includegraphics[height=6cm]
 {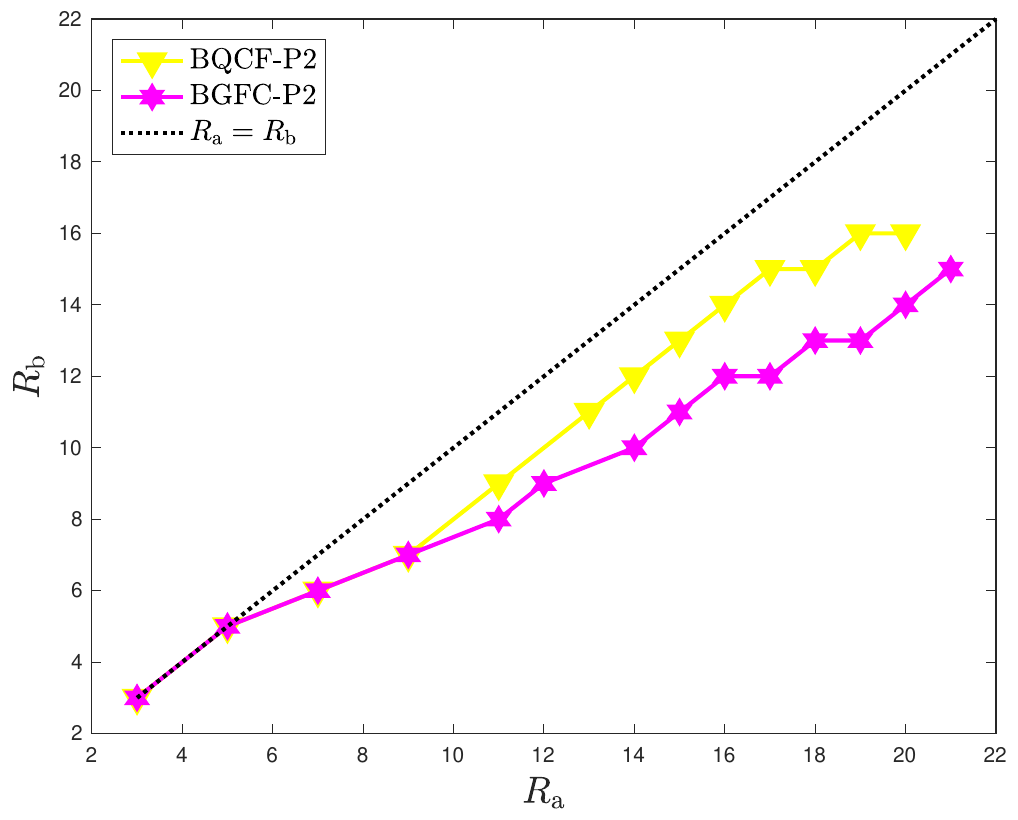}
	\caption{The relationship between the radius of the atomistic region $R_{\a}$ and the width of the blending region $R_{\b}$ in the adaptive computations for Frenkel defect.}
	\label{figs:RaRb_frenkel}
\end{center}
\end{figure}

We plot the relationship between $R_{\a}$ and $R_{\b}$ in the adaptive computations for Frenkel defect in Figure~\ref{figs:RaRb_frenkel}, where $R_{\a}$ is the radius of the atomistic region while $R_{\b}$ represent the width of the blending region. Our adaptive algorithm can achieve the optimal relationship between $R_{\a}$ and $R_{\b}$ for the case of Frenkel defect, indicating the robustness of main adaptive algorithm.

\begin{figure}[htb]
\begin{center}
	\includegraphics[height=6cm]{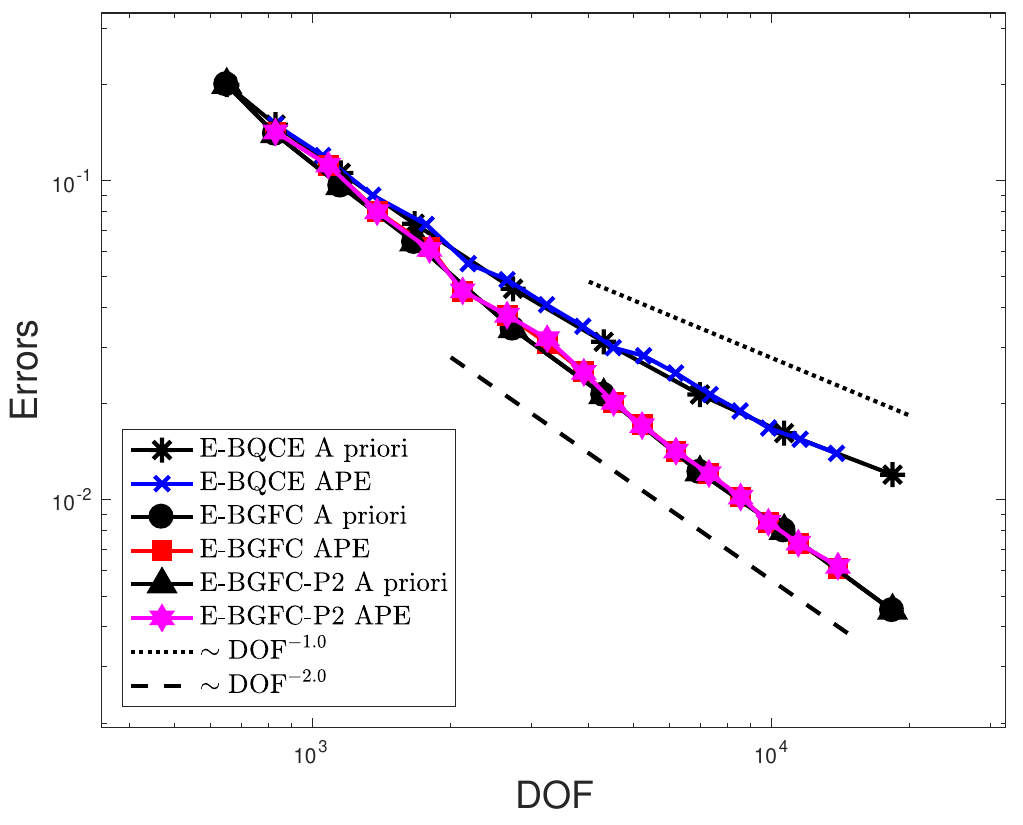}
	\caption{The convergences of energy error for various a/c coupling methods with respect to the number of degrees of freedom for Frenkel defect.}
	\label{fig:conv_E_frenkel_conv}
\end{center}
\end{figure}

As for the energy error, we observe in Figure~\ref{fig:conv_E_frenkel_conv} that all adaptive computations achieve the same optimal convergence rate ($N^{-1.0}$ for BQCE whereas $N^{-2.0}$ for BGFC-P1 and BGFC-P2) compared with the {\it a priori} graded mesh given by \cite{colz2016, fang2020blended}, which indicate the error estimator for energy error is still reliable.

\section{Conclusion and Outlook}
\label{sec: conclusion}

In this work we propose a unified framework of the residual based {\it a posteriori} error estimates and design the corresponding adaptive algorithms which is in principle suitable for general consistent multiscale coupling methods. We prove that the error estimator based on the residual forces can provide the upper bound of the true approximation error. As prototypical examples, we present a range of adaptive computations based on this reliable error estimator for the blended atomistic-to-continuum (a/c) coupling methods including the energy-based blended quasi-continuum (BQCE), the force-based blended quasi-continuum (BQCF) and the recently developed blended ghost force correction (BGFC) methods. We develop some coarse-grained techniques for the efficient evaluation of the error estimator and test with different types of crystalline defects some of which are not considered in previous related literature of the adaptive a/c coupling methods. The various numerical results demonstrate that, compared with the {\it a priori} error estimate, the adaptive algorithm leads to the same optimal convergence rate of the error with considerable computational efficiency. We point out that the techniques and strategies developed in this work do not require too much specific choice of the multiscale coupling scheme as long as the method is consistent.

Although we believe that the adaptive algorithm proposed in this paper is generally applicable for other common multiscale coupling schemes and more complex crystalline defects, this research still raise a few open problems which deserve further mathematical analysis and algorithmic developments.

\begin{itemize}    
    \item {\it More complex crystalline defects}: The main bottleneck for more complex defects such as dislocations and cracks is the construction and implementation of the corresponding a/c coupling methods. It is very difficult (if not impossible) to give rigorous {\it a priori} analysis for such problems. However, it should be the place where the adaptive a/c coupling methods reveal its most advantage and power and we believe that the current research makes it possible to develop efficient and robust adaptive algorithms for those problems with practical interests. 

    \item {\it Three dimensional model problems}: We expect to extend the current two dimensional setting to three dimensions. The adaptive algorithm proposed in this work has no fundamental obstacle but the implementation needs additional substantial effort since the highly non-trivial 3D mesh generation and adaptation are required. The recent advance \cite{fu2021adaptive} should provides good reference in this direction.  

    \item {\it More robust and flexible refinement strategy}: Although most of the steps in our analysis and algorithm are independent of the geometry of the material defect and computational domain, the adaptive algorithm can only adjust the radius of the atomistic region to refine the model. This limitation prevents us from capturing significant anisotropy in the defect core, elastic field, or defect nucleation. To address this issue and consider such generalizations, we need to evolve the interface anisotropically. One approach is to treat this as a free interface problem based on the error distribution, which may lead to a more robust implementation of model adaptivity. Recent advancements in the literature of QM/MM coupling methods provide one such example~\cite{adapqmmm2022}.
    
\end{itemize}

Both the theoretical and practical aspects discussed above will be explored in future work.

\appendix

\def\uH{u_h}

\section{Derivation and connection of atomistic stress}
\label{sec:app:Astress}

In this section, motivated by the stress formulation of the a/c coupling methods and corresponding stress based {\it a posteriori} error estimators~\cite{liao2018posteriori, wang2018posteriori}, we will make the connection between the mechanical notion of stress and the {\it a posteriori} estimator defined through $\phi$ in \eqref{eq:defnphi}. 

The atomistic stress formulation plays a significant role in the analysis of a/c coupling methods. Given displacement $u\in\UsH$, the second Piola stress tensor $\sigma(u)$ of atomistic model satisfies
\[
\<\delta\E^{\a}(u),v\> = \int_{\R^2} \sigma(u) \cdot \nabla v \dx, \quad \forall v\in\UsH.
\]

Next, we derive the concrete formulation of $\sigma(u)$ in the following. The first variation of the atomistic energy functional \eqref{energy-difference} is given by
\begin{eqnarray}\label{eq:var_QM_}
\< \delta\E^{\a}(u), v \>=\sum_{\ell \in \L} \sum_{\rho \in \Rg} V_{\ell, \rho}\big(Du(\ell)\big) \cdot D_{\rho}v(\ell), \quad \forall v\in\UsH. 
\end{eqnarray}
To derive its continuous version from which $\sigma(u)$ is then obtained, we apply the so-called localization formulation introduced in~\cite{OrtnerTheil2012}, that is, 
\begin{align}\label{eq:localization}
D_{\rho} \tilde{v}(\ell)=\int_0^1\nabla_{\rho}\tilde{v}(\ell+t\rho)\dt
=\int_{\R^{2}}\int_0^1\zeta(\ell+t\rho-x)\cdot \nabla_{\rho}v(x)\dt\dx,
\end{align}
where $\tilde{v}$ is a (quasi-)interpolation of $v$. We omit the construction here for the brevity of presentation and refer to~\cite{OrtnerTheil2012, AlexIntp:2012} for its detail.

Due to the discreteness nature of atomistic models as well as the proposition of $\tilde{v}$~\cite{AlexIntp:2012}, we can evaluate \eqref{eq:var_QM_} by replacing the test function $v$ with $\tilde{v}$, which then leads to
\begin{align}\label{eq:var_QM}
\< \delta\E^{\a}(u), \tilde{v} \>&=\sum_{\ell \in \L} \sum_{\rho \in \Rg} V_{\ell, \rho}\big(Du(\ell)\big) \cdot D_{\rho}\tilde{v}(\ell) \nonumber \\
&=\int_{\R^2} \sum_{\ell \in \L} \sum_{\rho \in \Rg} \rho \otimes V_{\ell, \rho}\big(Du(\ell)\big) \int_0^1\zeta(\ell+t\rho-x)\dt \cdot \nabla v(x) \dx   \nonumber \\
&=: \int_{\R^2}\sigma(u)(x) \cdot \nabla v(x) \dx.
\end{align}
Hence, the concrete atomistic stress reads
\begin{eqnarray}\label{eq:QM_stress}
\sigma(u)(x) := \sum_{\ell \in \L} \sum_{\rho \in \Rg} \rho \otimes V_{\ell, \rho}\big(Du(\ell)\big)\cdot \chi_{\ell,\rho}(x)
\end{eqnarray}
with the ``smeared bonds'' $\chi_{\ell,\rho}(x) := \int_0^1 \zeta(\ell+t\rho-x)\dt$. Recall that $u_h$ is the solution of blended a/c coupling methods, we have $I_{\a}u_h\in \UsH$. Then, the {\it a posteriori} atomistic stress tensor $\sigma(I_{\a}u_h)$ is acquired. However, it can not be applied in practical adaptive algorithm since $\sigma(I_{\a}u_h)(x)$ is relatively difficult to be evaluated through the ``smeared bonds'' $\chi_{\ell, \rho}(x)$.

Moreover, we could rewrite the first variation of atomistic energy functional by considering the residual forces (cf.~\eqref{eq:forcenew}), namely
\begin{align}\label{eq:var_QMf}
\< \delta\E^{\a}(I_{\a}\uH), \tilde{v} \>
=& \sum_{\ell \in \L} \sum_{\rho \in \Rg} \Big( V_{\ell-\rho, \rho}\big(D\uH(\ell-\rho)\big) - V_{\ell, \rho}\big(D\uH(\ell)\big) \Big) \cdot \tilde{v}(\ell) \nonumber \\
=& \int_{\R^{2}} \widehat{\F}^{\a}(I_{\a}\uH)(x) \cdot v(x) \dx.
\end{align}

Combining \eqref{eq:var_QMf} with \eqref{eq:var_QM}, one can obtain the following equation
\begin{eqnarray}\label{eq:stress_force}
\int_{\R^2}\sigma(I_{\a}\uH)(x) \cdot \nabla v(x)\dx = \int_{\R^2} \widehat{\F}^{\a}(I_{\a}\uH)(x) \cdot v(x)\dx, \quad \forall v \in \UsH.
\end{eqnarray}
We note that it is highly related to the equation from which the {\it a posteriori} error estimator $\|\nabla\phi\|_{L^2(R^2)}$ is obtained. As a matter of fact, according to the Helmholtz-Hodge decomposition \cite{paula2016, Or:2011a}, $\sigma(I_{\a}\uH)$ can be decomposed as a sum of two orthogonal components:
\begin{equation}
\label{eq:HHD}
\sigma(I_{\a}\uH) = \nabla \phi + \nabla\times \psi,
\end{equation}
with $\phi \in \UsH$, $\psi\in \UsH$.
$\nabla\phi\in L^2$ is called the ``curl-free" component, and $\nabla\times \psi \in  L^2$ is divergence-free in the weak sense, i.e., $\int_{\R^2}\nabla\times\psi(x) \cdot \nabla v(x) \dx = 0$. 

Combining Theorem \ref{th:mainresult} and Lemma \ref{th:ctsphi}, we have
\begin{align*}
\eta^2(\uH) \lesssim \|\mR(I_{\a}u_h)\|^2_{(\UsH)^{*}} \lesssim \|\nabla \phi\|_{L^2}^2\leq \|\nabla\phi\|_{L^2}^2 + \|\nabla\times\psi\|_{L^2}^2=\|\sigma(I_{\a}\uH)\|_{L^2}^2.
\end{align*}
Therefore, $\|\sigma(I_{\a}\uH])\|_{L^2}$ provides an upper bound for the approximation error, which is consistent with the result shown in~\cite[Theorem 3.7]{wang2018analysis}.

In \eqref{eq:HHD}, we can uniquely define $\phi$ by $\Delta \phi = \nabla \cdot \sigma(I_{\a} \uH)$ (in the weak sense). On the other hand, we can choose an arbitrary divergence-free component $\nabla\times\psi$ in $\sigma(I_{\a} \uH)$ to satisfy \eqref{eq:stress_force}. Therefore, $\sigma(I_{\a} \uH)$ in the sense of \eqref{eq:stress_force} is not unique. Motivated by the stress tensor correction~\cite[Eq.(70)]{wang2018analysis}, we consider the following problem to obtain a uniquely defined QM stress tensor,
\begin{equation}\label{eq:minSigma}
\bar{\psi} \in \arg\min_{\psi \in  \UsH} \big\{ \|\sigma(I_{\a}\uH)\|_{L^2} = \|\nabla\phi + \nabla\times\psi\|_{L^2} \big\}.
\end{equation}

A straightforward calculation and the orthogonality of the two components of Helmholtz-Hodge decomposition lead to $\nabla\times\bar{\psi} = 0$. Hence, we denote the corresponding uniquely defined QM stress tensor as
\begin{equation}\label{eq:stress_grad}
\sigma^0(I_{\a} \uH) := \nabla\phi, \quad {\rm for}\; \phi \in \UsH.
\end{equation}

As a result, we recover the equation \eqref{eq:defnphi} which is used in the a posteriori estimates in Section~\ref{sec:sub:sub:imp_stress}.

\section{Proof of Theorem~\ref{th:mainresult}}
\label{sec:app:analysis}

First of all, we introduce the definition of truncation operator. To estimate the approximation error introduced by this discretization we particularly need to account for the truncation of the domain. We assume for the sake of technical convenience, that there exists a radius $\RO$ such that $B_{\RO} \subset\Omega\subset B_{2\RO}$; that is, $\Omega$ is approximately  a ball. This allows us to define a simple truncation operator, following \cite{2013-defects}, $T_{\RO}: \dot{H}^1(\R^d)\rightarrow H^1_0(\Omega)$, 
\begin{align}\label{eq:trmod}
T_{\RO}v(x) := \eta(x)\big(v(x)-a_{\RO}\big),
\qquad a_{\RO} = \mint_{B_{\RO}\setminus B_{\RO/2}}v(x) \dx
\end{align}
where $\eta$ is a $C^1$ cut-off function; $\eta = 1$ in $B_{\RO/2}$, $\eta = 0$ in $B_{\RO}^{\rm c}$ and $|\nabla \eta| \lesssim \RO^{-1}$. Following \cite{2013-defects}, for $v_{\RO} = T_{\RO} v$ and $a_{\RO}$ defined by \eqref{eq:trmod} we readily obtain the estimates
    \begin{align}\label{eq:trwgw}
	\| v - v_{\RO} - a_{\RO} \|_{L^2(\Omega)} &\lesssim \RO \| \nabla v - \nabla v_{\RO} \|_{L^2(\Omega \setminus B_{\RO})}, \qquad \text{and} \\
\label{eq:trwgw2}
	\| \nabla v - \nabla v_{\RO} \|_{L^2(\Omega)} &\lesssim  \| \nabla v \|_{L^2(\Omega \setminus B_{\RO/2})}.
   \end{align}

\begin{proof}
First of all, we estimate the term $\| \nabla \phi - \nabla \phi_{\rm a} \|_{L^2(\R^2)}$. By \eqref{eq:defnphi}, for $v\in\dot{H}^1$, we have
\begin{align}\label{eq:phi-phia}
		\| \nabla \phi - \nabla \phi_{\rm a} \|_{L^2(\R^2)}
		 &= \sup_{\| \nabla v \|_{L^2} = 1}
\int_{\R^2} \big( \nabla \phi \cdot \nabla v - \nabla \phi_{\a} \cdot \nabla v \big) \dx \nonumber \\
		&= \sup_{\| \nabla v \|_{L^2} = 1}
\int_{\R^2} \big( \widehat{\F}^{\a} \cdot v - \nabla \phi_{\a} \cdot \nabla v \big) \dx.
\end{align}
Let $v_R := T_{R_\Omega} v$ be defined by \eqref{eq:trmod}, then we split the residual into two parts, 
\begin{align}
	\int_{\R^2} \big( \widehat{\F}^{\a} \cdot v - \nabla \phi_{\rm a} \cdot \nabla v  \big) \dx
		&= \int_{\R^2} \widehat{\F}^{\a} \cdot (v-v_R) \dx  - \int_{\R^2} \nabla \phi_{\rm a} \cdot (\nabla v -\nabla v_R)\dx \nonumber \\
		&=: T_1 + T_2. 
\end{align}

To estimate the term $T_1$, we have
\begin{align}\label{eq:EstT1}
T_1 = \int_{\R^2} \widehat{\F}^{\a} \cdot (v-v_R) \dx \lesssim~R_\Omega \cdot \|\widehat{\F}^{\a}\|_{L^2(\Omega\setminus B_{R_\Omega/2})} \cdot \|\nabla v\|_{L^2}.
\end{align}
where the last inequality follows from the estimate \eqref{eq:trwgw2}. 

Similarly, for the term $T_2$, we can obtain
\begin{align}\label{eq:EstT2}
T_2 = \int_{\R^2} \nabla \phi_{\rm a} \cdot (\nabla v -\nabla v_R)\dx \lesssim \|\nabla \phi_\a\|_{L^2(\Omega\setminus B_{R_\Omega/2})}\cdot \|\nabla v\|_{L^2}.
\end{align}

Hence, combining \eqref{eq:phi-phia}, \eqref{eq:EstT1} and \eqref{eq:EstT2}, one can acquire 
\begin{align}\label{eq:Est-phi-phia}
		\| \nabla \phi - \nabla \phi_{\rm a} \|_{L^2(\R^2)}
		 \lesssim R_\Omega \cdot \|\widehat{\F}^{\a}\|_{L^2(\Omega\setminus B_{R_\Omega/2})} + \|\nabla \phi_\a\|_{L^2(\Omega\setminus B_{R_\Omega/2})} =: \rho_{\rm tr}(u_h).
\end{align}

According to Lemma \ref{lemma:res-F},
we already know that $ \| {\textsf{R}}(I_{\rm a} u_{h}) \|_{(\UsH)^*} \simeq \| \nabla \phi \|_{L^2}$. Moreover, we can deduce that $\| \nabla \phi \|_{L^2} =  \| \widehat{\F}^{\a} \|_{(\dot{H}^1)^*}$. Since $\phi$ and $\phi_\a$ are the solutions of \eqref{eq:defnphi} and \eqref{eq:stress_force_trun_T}, respectively, we use Galerkin orthogonality to write
\begin{align}
    \label{eq:phi-phih}
    \| \nabla \phi \|^2_{L^2} - \| \nabla \phi_{\a}\|^2_{L^2}  =& \<\nabla \phi - \nabla \phi_{\a}, \nabla \phi - \nabla \phi_{\a}\> + 2\<\nabla \phi - \nabla \phi_{\a}, \nabla \phi_\a\> \nonumber \\
    =& \| \nabla \phi - \nabla \phi_{\a}\|^2_{L^2}.
\end{align}

Hence, we have $\| \nabla \phi \|^2_{L^2} \geq \| \nabla \phi_{\a}\|^2_{L^2}$, which immediately indicates that $\eta(u_h) \lesssim \| \mR(I_{\rm a} u_{h}) \|_{(\UsH)^*}$.

To obtain an upper bound for $\| \nabla \phi \|_{L^2}$ we use Cauchy's inequality to estimate
\begin{align*}
	\| \nabla \phi \|^2_{L^2}
	=
	\| \nabla \phi_{\a}\|^2_{L^2}
	+
	\| \nabla \phi - \nabla \phi_{\a}\|^2_{L^2}
	\leq
	\big( \| \nabla \phi_{\a}\|_{L^2}
	+
	\| \nabla \phi - \nabla \phi_{\a}\|_{L^2} \big)^2.
\end{align*}
We deduce $\| \mR(I_{\rm a} u_{h}) \|_{(\UsH)^*} \lesssim \eta(u_h) + \rho_{\rm tr}(u_h)$, which finishes the stated result.
\end{proof}

\bibliographystyle{plain}
\bibliography{qc}
\end{document}

%% file: notation.tex
\renewcommand{\cases}[1]{\left\{ \begin{array}{rl} #1 \end{array} \right.}
\newcommand{\smfrac}[2]{{\textstyle \frac{#1}{#2}}}
\newcommand{\myvec}[1]{\left[ \begin{vector} #1 \end{vector} \right]}
\newcommand{\mymat}[1]{\left[ \begin{matrix} #1 \end{matrix} \right]}

\def\Xint#1{\mathchoice
{\XXint\displaystyle\textstyle{#1}}%
{\XXint\textstyle\scriptstyle{#1}}%
{\XXint\scriptstyle\scriptscriptstyle{#1}}%
{\XXint\scriptscriptstyle\scriptscriptstyle{#1}}%
\!\int}
\def\XXint#1#2#3{{\setbox0=\hbox{$#1{#2#3}{\int}$ }
\vcenter{\hbox{$#2#3$ }}\kern-.6\wd0}}
\def\mint{\Xint-}

\def\b{\big}
\def\B{\Big}
\def\bg{\bigg}
\def\Bg{\Bigg}

\def\sep{\,|\,}
\def\bsep{\,\b|\,}
\def\Bsep{\,\B|\,}

\def\diam{{\textrm{diam}}}
\def\conv{{\textrm{conv}}}
\def\t{\top} 
\def\sign{{\textrm{sgn}}}
\def\id{{\textrm{id}}}
\def\supp{{\textrm{supp}}}
\def\diam{{\textrm{diam}}}

\def\R{\mathbb{R}}
\def\N{\mathbb{N}}
\def\Z{\mathbb{Z}}
\def\C{\mathbb{C}}
\def\bbV{\mathbb{V}}

\def\WW{W}
\def\CC{C}
\def\HH{H}
\def\LL{L}
\def\DD{{D}'}
\def\Ys{\mathscr{Y}}

\def\WWh{\dot{W}}
\def\Ycb{Y}
\def\WWhz{\dot{W}_0}
\def\Ycbz{\Ycb_0}

\def\dx{{\,{\textrm{d}}x}}
\def\dy{\,{\textrm{d}}y}
\def\dz{\,{\textrm{d}}z}
\def\dr{\,{\textrm{d}}r}
\def\dt{\,{\textrm{d}}t}
\def\ds{\,{\textrm{d}}s}
\def\dd{\textrm{d}}
\def\pp{\partial}
\def\dV{\,\textrm{dV}}
\def\dA{\,{\textrm{dA}}}
\def\db{\,{\textrm{db}}}
\def\dlam{\,{\textrm{d}}\lambda}

\def\<{\langle}
\def\>{\rangle}

\def\ol{\overline}
\def\ul{\underline}
\def\ot{\widetilde}
\newcommand{\ut}[1]{\underset{\widetilde{\hspace{2.5mm}}}{#1}}

\def\mA{{\textsf{A}}}
\def\mB{\textsf{B}}
\def\mC{\textsf{C}}
\def\mF{\textsf{F}}
\def\mG{\textsf{G}}
\def\mH{\textsf{H}}
\def\mI{\textsf{I}}
\def\mJ{\textsf{J}}
\def\mP{\textsf{P}}
\def\mQ{\textsf{Q}}
\def\mR{\textsf{R}}
\def\mM{\textsf{M}}
\def\mS{\textsf{S}}
\def\mO{\textsf{0}}
\def\mL{\textsf{L}}

\def\sym{\textsf{sym}}
\def\tr{\textsf{tr}}
\def\el{\textsf{el}}

\def\bfa{\textbm{a}}
\def\bfg{\textbm{g}}
\def\bfrho{\mathbm{\rho}}
\def\bfv{\textbm{v}}
\def\bfh{\textbm{h}}
\def\bfO{\textbm{0}}

\def\bbA{\mathbb{A}}
\def\bbB{\mathbb{B}}
\def\bbC{\mathbb{C}}
\def\bbI{\mathbb{I}}

\def\Hs{\mathcal{H}}


\newcommand{\transpose}{{\!\top}}
\newcommand{\Da}[1]{D_{\!#1}}
\newcommand{\Dc}[1]{\D_{#1}}
\def\D{\nabla}
\def\del{\delta}
\def\ddel{\delta^2}
\def\dddel{\delta^3}

\def\loc{\textrm{loc}}

\def\qc{\textrm{qc}}
\def\c{\textrm{c}}
\def\h{\textrm{h}}

\def\eps{\varepsilon}
\def\tot{\textrm{tot}}
\def\cb{\textrm{cb}}
\def\a{\textrm{a}}
\def\c{\textrm{c}}
\def\ac{\textrm{ ac}}
\def\i{\textrm{i}}
\def\nn{\textrm{nn}}
\def\refl{\textrm{rfl}}
\def\qnl{\textrm{qnl}}
\def\stab{\textrm{stab}}
\def\conv{\textrm{conv}}
\def\supp{\textrm{supp}}

\def\L{\Lambda}
\def\Is{\mathcal{I}}
\def\oIs{\ol{\Is}}
\def\As{\mathcal{A}}
\def\Cs{\mathcal{C}}
\def\Fs{\mathcal{F}}
\def\Ks{\mathcal{K}}
\def\Us{\mathscr{U}}
\def\Usz{\Us_0}
\def\Ush{\dot{\Us}^{1,2}}
\def\Ushd{\dot{\Us}^{-1,2}}
\def\Usp{{\Us}^{1,p}}
\def\Usc{\Us^c}

\def\Bs{\mathcal{B}}
\def\Ls{\mathcal{L}}
\def\bbL{\mathbb{L}}

\def\yF{y_\mF}
\def\uF{u_\mF}

\def\E{\mathcal{E}}
\def\Ea{\E^\a}
\def\Eb{\E^\textrm{b}}
\def\Ha{H^\a}
\def\Ei{\E^\i}
\def\Ec{\E^\c}
\def\Eh{\E_\h}
\def\F{\mathscr{F}}
\def\Hc{H^\c}
\def\Eqnl{\E^\qnl}
\def\Hqnl{H^\qnl}
\def\Erefl{\E^\refl}
\def\Hrefl{H^\refl}
\def\Eac{\E^\ac}
\def\Hac{H^\ac}
\def\Estab{\E^\stab}
\def\Hstab{H^\stab}

\def\dW{W'}
\def\ddW{W''}

\def\RO{\mathcal{R}}

\def\Es{\Phi}

\def\Eatot{\E^\a_\tot}
\def\Eatot{\E^\c_\tot}

\def\Om{{\R^d}}
\def\Vi{V^\i}
\def\Vc{V^\c}
\def\Vs{\mathscr{V}}

\def\tily{\tilde y}
\def\tilz{\tilde z}
\def\tilu{\tilde u}
\def\tilv{\tilde v}
\def\tile{\tilde e}
\def\tilw{\tilde w}
\def\tilf{\tilde f}

\def\bary{y}
\def\barz{z}
\def\barv{v}
\def\barw{w}
\def\baru{u}
\def\bare{e}
\def\barf{f}


\def\ve{\varepsilon}
\def\L{\Lambda}
\def\La{\L^{\a}}
\def\Li{\L^{\i}}
\def\Lc{\L^{\c}}
\def\yd{y_0}
\def\Nhd{\mathcal{N}}
\def\rcut{r_{\textrm{cut}}}
\def\Rg{\mathscr{R}}
\def\Rgp{\Rg^{+}}
\def\vsig{\varsigma}
\def\T{\mathcal{T}}
\def\Tp{T^+}
\def\Tm{T^-}
\def\np{\nu^+}
\def\nm{\nu^-}
\def\T{\mathcal{T}}
\def\Th{\mathcal{T}_\h}
\def\Ta{\T_\a}
\def\Te{\mathscr{T}_\varepsilon}
\def\Fc{\mathscr{F}_\c}
\def\Fi{\mathscr{F}_\i}
\def\Fh{\mathscr{F}_\h}
\def\UsT{\Us_\h}
\def\Ih{I_\h}
\def\Ia{I_\a}
\def\Ie{I_{\varepsilon}}
\def\Tmu{\mathscr{T}_\mu}
\def\vor\textrm{vor}
\def\s{\sigma}
\def\sh{\sigma^\h}
\def\sa{\sigma^\a}
\def\sc{\sigma^\c}
\def\Oma{\Omega^\a}
\def\PO{\textrm{P}_0}
